\documentclass[11pt]{article}

\usepackage{makeidx}
\usepackage{amssymb}
\usepackage{amsfonts}
\usepackage{amsmath}
\usepackage{euscript}
\usepackage{theorem}
\usepackage{enumerate}

\theorembodyfont{\rm}      

\theoremstyle{change}      

\begingroup

\newtheorem{thm}{Theorem\hskip 5mm}[section]

\newtheorem{theorem}{Theorem\hskip 5mm}[section]

\newtheorem{prop}[thm]{Proposition\hskip 5mm}

\newtheorem{lem}[thm]{Lemma\hskip 5mm}

\newtheorem{notn}[thm]{Notation\hskip 5mm}

\newtheorem{defn}[thm]{Definition\hskip 5mm}

\newtheorem{note}[thm]{Note\hskip 5mm}

\endgroup

\begingroup

\endgroup

\begingroup

\endgroup

\def\s{\sigma}
\def\vf{{\varphi}}
\def\GL{\mathrm{GL}}
\def\Sp{\mathrm{Sp}}
\def\O{\mathrm{O}}

\def\B{{\mathcal B}}
\def\Bil{\mathrm{Bil}}
\let\la=\langle  \let\ra=\rangle

\def\even{{\rm even}}
\def\odd{{\rm odd}}
\def\tdeg{{\rm tdeg}}
\def\ndeg{{\rm ndeg}}
\def\GI{{G[V_\infty]}}

\def\GIK{{G[V_\infty]\cap G[{}^\infty V/V_\infty]}}
\def\KB{{G[V_\infty]\cap G[{}^\infty V/V_\infty]/G[{}^\infty V]}}
\def\BB{{G[{}^\infty V]}}
\def\R{{\rm Rad}}

\def\VI{{{}^\infty V}}

\def\WI{{V/{}^\infty V}}

\def\t{{\bf t}}
\def\l{{\mathcal L}}
\def\m{{\mathcal M}}
\def\P{{\mathcal P}}
\def\dim{\mathrm {dim}\,}
\def\End{\mathrm {End}}
\def\Hom{\mathrm {Hom}}
\def\ker{\mathrm {ker}\,}

\def\so{{\mathfrak s\mathfrak o}}
\def\o{{\mathfrak o}}
\def\sp{\mathfrak{ sp}}

\def\DJ{{\hbox{D\kern-.8em\raise.15ex\hbox{--}\kern.35em}}}
\def\DJo{\DJ okovi\'c}

\begin{document}

\begin{center}
{\bf STRUCTURE OF THE GROUP PRESERVING A BILINEAR FORM}
\end{center}
\begin{center}
{\sc Fernando Szechtman}
\end{center}
\begin{center}
{\small\textit{Department of Mathematics \& Statistics, University
of Regina, Saskatchewan, Canada, S4S 0A2\\ e-mail:
szechtf@math.uregina.ca }}
\end{center}

\begin{abstract}
\noindent We study the group of all linear automorphisms
preserving an arbitrary bilinear form.
\end{abstract}

\section{Introduction}

Let $\vf:V\times V\to F$ be a bilinear form on a finite
dimensional vector space $V$ over a field $F$. We refer to the
pair $(V,\vf)$ as a bilinear space. The goal of this paper is to
describe the structure of the group $G=G(V)$ of all linear
automorphisms of $V$ preserving $\varphi$.

Classical groups such as the general linear, symplectic and
orthogonal groups arise in this fashion. These classical cases
have been thoroughly investigated (see \cite{JD2}, \cite{HO}).
Arbitrary non-degenerate forms possess an asymmetry (as defined in
\cite{CR}), which exerts a considerable influence on the structure
of the bilinear space and associated group. This has recently been
exploited by J. Fulman and R. Guralnick \cite{FG}, where an array
of useful information about $G$ is presented. The study of $G$ for
a general bilinear form, possibly degenerate, over an arbitrary
field does not seem to have been hitherto considered. The presence
of a degenerate part enriches the structure of $G$, and it is in
this regard that our main contribution takes place.

In general terms our approach consists of extracting structural
information about $G$ by examining how $G$ acts on $V$ and its
various $FG$-submodules.

Knowledge of the structure of $V$ as an $FG$-module will therefore
be necessary. Our references in this regard will consist of the
paper \cite{CR} by C. Riehm, its appendix \cite{PG} by P. Gabriel,
and our recent article \cite{DS1} with D. Djokovic.

An important decomposition of $V$ to be considered is
$$
V=V_\odd\perp V_\even\perp V_\ndeg,
$$
where $V_\odd$, respectively $V_\even$, is the orthogonal direct
sum of indecomposable degenerate bilinear spaces of odd,
respectively even, dimension, and $V_\ndeg$ is non-degenerate. We
identify $G(V_\odd)$, $G(V_\even)$ and $G(V_\ndeg)$ with subgroups
of $G(V)$ by means of this decomposition.

As noted in \cite{DS1}, while  $V_\odd$, $V_\even$ and $V_\ndeg$
are uniquely determined by $V$ up to equivalence of bilinear
spaces, they are not unique as subspaces of $V$, and in particular
they are not $G$-invariant. Thus, one attempt to understand $G$
would consist of studying the structure of $G(V_\odd)$,
$G(V_\even)$ and $G(V_\ndeg)$ separately, and then see how these
groups fit together to form $G$.

This approach turns out to be fruitful, as we proceed to describe,
with the notable exception of the structure of $G(V_\ndeg)$ in the
special case when the asymmetry of $V_\ndeg$ is unipotent and the
underlying field $F$ has characteristic 2. This is what C. Riehm
refers to as Case IIb in his paper.

We begin our journey in section 2 by establishing notation and
terminology. Two important $G$-invariant subspaces of $V$
described here are what we denote by $V_\infty$ and ${}^\infty V$
in \cite{DS1}. Briefly, any choice for $V_\odd$ will contain
$V_\infty$, which is in fact the only totally isotropic subspace
of $V_\odd$ of maximum possible dimension; whatever the choices
for $V_\even$ and $V_\ndeg$, it turns out that ${}^\infty
V=V_\infty\perp V_\even\perp V_\ndeg$.

Section 3 contains basic tools regarding $V$ and $G$ to be used
throughout the paper. An important feature of this section is the
introduction -in Definition \ref{bert1}- of a family of
1-parameter subgroups of $G(V_\odd)$ which will play a decisive
role in shedding light on the structure of both $G(V_\odd)$ and
$G$.

The actual paper can be said to begin in section 4. We first
introduce a normal subgroup $N$ of $G$, defined as the
intersection of various pointwise stabilizers in $G$ when it acts
on  certain sections of the $FG$-module $V$. It is shown in
Theorem \ref{E} that $G=N\rtimes E$, where $E\cong \underset{1\leq
i\leq t}\Pi \GL_{m_i}(F)$ and the parameters $t$ and $m_1,...,m_t$
depend only on $V$, as explained below. We have
$$
V_{\odd}=V_1\perp V_2\perp\cdots\perp V_t,
$$
where each $V_i$ is the orthogonal direct sum of $m_i$ bilinear
subspaces, each of which is isomorphic to a  Gabriel block of size
$2s_i+1$, with $s_1>s_2>\cdots>s_t$. Here by a Gabriel block of
size $r\geq 1$ we mean the only indecomposable degenerate bilinear
space of dimension $r$, up to equivalence, namely one that admits
a nilpotent Jordan block of size $r$ as its Gram matrix.

As a byproduct of the results in this section we are able to
describe the irreducible constituents of the $FG$-submodule
$V_\infty$ of $V$. This is taken up in section 5 (Theorem
\ref{irrvi}). These constituents are seen later in section 6 to be
intimately connected to certain $FG$-modules arising as sections
of $G$ itself. The end of this section also gives a second
decomposition for $G$, namely
$$
G(V)=G[{}^\infty V/V_\infty]\rtimes (G(V_\even)\times G(V_\ndeg)),
$$
where in general $G[Y]$ denotes the pointwise stabilizer of $G$
acting on a $G$-set $Y$.

Section 6 goes much deeper than previous sections. In view of the
decomposition $G=N\rtimes E$ and the clear structure of $E$, our
next goal is to study $N$ and the action of $G$ upon it. In
Theorem \ref{U} we prove $N=G[V_\infty]\rtimes U$, where $U$ is a
unipotent subgroup of $G(V_\odd)$ generated by the 1-parameter
subgroups referred to above. An extensive analysis of the
nilpotent group $N/G[V_\infty]$ is carried out. First of all, its
nilpotency class is seen in Theorem \ref{nilpdegree} to be $t-1$.
As a nilpotent group $N/G[V_\infty]$, possesses a descending
central series. We actually produce in Theorem \ref{main} a
$G$-invariant descending central series for $N/G[V_\infty]$ each
of whose factors has a natural structure of $FG$-module, and
irreducible at that. These irreducible $FG$-modules are in close
relationship with the irreducible constituents of the $FG$-module
$V_\infty$. By taking into account all factors in our series we
deduce a formula for the dimension of $U$, in the
algebraic/geometric sense, which turns out to be equal to the
number of 1-parameter groups generating $U$ and referred to above.
Theorem \ref{dim-U} proves
$$\dim U=\underset{1\leq i<j\leq t}\sum (s_i-s_j+1) m_im_j.$$

As a byproduct of the results in this section we also obtain in
Theorem \ref {irrvI} the irreducible constituents of the
$FG$-module $V/{}^\infty V$, which in turn are closely related to
other $FG$-modules also arising as sections of $G$.

Section 7 concentrates on the next logical target, namely
$G[V_\infty]$. We know from above that $G=N\rtimes E$ and
$N=G[V_\infty]\rtimes U$. Here we prove (Theorem \ref{bigsplit})
that $U$ actually normalizes $E$ -so $G=G[V_\infty]\rtimes
(U\rtimes E)$-, that $G[V_\infty]$ admits the decomposition
$G[V_\infty]=(G[V_\infty]\cap G[{}^\infty V/V_\infty])\rtimes
(G(V_{\even})\times G(V_{\ndeg}))$, and that $U\rtimes E$ actually
commutes with $G(V_{\even})\times G(V_{\ndeg})$ elementwise. We
thus obtain the important decomposition
$$
G=(G[V_\infty]\cap G[{}^\infty V/V_\infty])\rtimes
\left(G(V_{\even})\times G(V_{\ndeg})\times (U\rtimes E)\right).
$$

With the structure of $U\rtimes E$ already clarified, the next
step consists of studying $G(V_{\even})$ and $G(V_{\ndeg})$ on
their own, and see what is the structure of $G[V_\infty]\cap
G[{}^\infty V/V_\infty]$.

Note that $G(V_\odd)$ seems to be absent above. But that is only
an illusion, which is clarified in section 11. In fact,
$G(V_\odd)$ is essentially what is holding the above decomposition
together. If $V_\odd=(0)$ the $V_\even$ and $V_\ndeg$ are in fact
$G$-invariant, as \cite{DS1} shows, so $G=G(V_\even)\times
G(V_\ndeg)$.

Section 8 begins by laying the foundations (Theorem \ref{YZ}) for
a combined attack on $G(V_{\even})$ and certain direct factors of
$G(V_{\ndeg})$. Theorem \ref{ee} then exploits this by describing
$G(V_\even)$ as the centralizer of a nilpotent element of known
similarity type in the general linear group.

Attention in section 9 is focused on $G(V_{\ndeg})$. This group is
approached via the study of $V_\ndeg$ as a module over the
polynomial algebra $F[t]$ by means of the asymmetry of
$\vf\vert_{V_\ndeg}$, as outlined in \cite{CR}. Thus (see equation
(\ref{dw})) $G(V_\ndeg)$ is isomorphic to the direct product of
groups of the form $G(W)$, where $W$ is a non-degenerate bilinear
space whose type, according to C. Riehm, is either I, IIa or IIb.

In the first case $G(W)$ is seen (in Theorem \ref{Ha} via Theorem
\ref{YZ}) to be isomorphic to the centralizer in a general linear
group of a linear automorphism of known similarity type. If $F$ is
algebraically closed this linear automorphism can be replaced by a
nilpotent endomorphism.

Case IIa is more difficult. We find (Theorem \ref{gral}) $G(W)$ to
be equal to the centralizer in a symplectic or orthogonal group of
a particular linear endomorphism. If $F$ has characteristic not 2
this element is in the corresponding symplectic or orthogonal Lie
algebra. If in addition $F$ is algebraically closed we can ensure
(Theorems \ref{sv+1} and \ref{sv-1}) that this element is
nilpotent of known similarity class. These centralizers are
described in various places, e.g. in \cite{JJ,SS}.

As mentioned already above the case when the asymmetry of
$V_\ndeg$ is unipotent and $F$ has characteristic 2, i.e. Case
IIb, remains unsolved.

Section 10 concentrates on $\GIK$. One sees rather rapidly (Lemmas
\ref{step2}, \ref{step3} and \ref{step4}) that $\GIK$ is unipotent
of nilpotency of class $\leq 2$ having $\BB$ in its center, the
corresponding quotient group being abelian. Thus the study of
$\GIK$ is divided into that of $\BB$ and $\GIK/\BB$.

Well, $\BB$ is naturally an $FG$-module and, much as in section 6,
we find (Theorem \ref{final}) its irreducible constituents and
explain how they relate to those of $V/{}^\infty V$. As in section
6, this requires considerable amount of work. In particular, the
dimension of $\BB$ is found. We also compute (Proposition
\ref{ufa}) the dimension of the quotient group $\GIK/\BB$, thereby
obtaining (Theorem \ref{ulto}) a formula for the dimension of
$\GIK$, which reads
$$
\dim\GIK=\dim (V/V_\infty)\times (m_1+\cdots+m_t).
$$

Section 11 furnishes a few more decompositions for $G$ and
$G(V_\odd)$ (Theorems \ref{kop} and \ref{kop2}) and includes an
example (Theorem \ref{homo}) on the structure of $G(V_\odd)$ in a
special but interesting case. The structure of $G(V_\odd)$ is
fully revealed in this case.

Our last section makes some comments on an alternative approach to
the study of $G$.

A few words about the origin of this paper are in order. After our
joint work \cite{DS1} with D. Djokovic, we were excited about the
prospect of being able to attack the present problem. We worked
rather intensively together for quite some time in fruitful
collaboration. Each of us built his own version of the paper, and
at one point our methods and some of our goals became too far
apart for us to be able to amalgamate them into a single paper.
Even though we agreed to submit our versions separately, the
outcome of this project should be regarded as joint work.

\section{Generalities}

Let $F$ be a field. A {\em bilinear space} over $F$ is a pair
$(V,\vf)$, where $V$ is a finite dimensional $F$-vector space and
$\vf :V\times V\to F$ is a bilinear form. An {\em isometry} from a
bilinear space $(V_1,\vf_1)$ to a bilinear space $(V_2,\vf_2)$ is
a linear isomorphism $g:V_1\to V_2$ satisfying
$$
\vf_2(gv,gw)=\vf_1(v,w), \quad v,w \in V_1.
$$
Two bilinear spaces are {\em equivalent} if there exists an
isometry between them.
The {\em isometry group} of a bilinear space $(V,\vf)$ is the
group of all isometries from $(V,\vf)$ into itself.

We henceforth fix a bilinear space $(V,\vf)$. Its isometry group
will be denoted by $G(V,\vf)$, $G(V)$, $G(\vf)$, or simply by $G$.
Explicit reference to $\vf$ will be omitted when no confusion is
possible. We shall often write $\la v,w\ra$ instead of $\vf(v,w)$.

The space of all bilinear forms on $V$ will be denoted by
$\Bil(V)$. There is an action of $\GL(V)$ on $\Bil(V)$ given by
$$
(g\cdot\phi)(v,w)=\phi(g^{-1}v,g^{-1}w),\quad g\in \GL(V),
\phi\in\Bil(V), v,w\in V.
$$
Thus the isometry group of $(V,\vf)$ is the stabilizer of $\vf$
under this action.

If $U$ is a subspace of $V$, then $U$ becomes a bilinear space by
restricting $\vf$ to $U\times U$. We write $V=U\perp W$ if
$V=U\oplus W$ and $\langle U,W\rangle=\langle W,U\rangle=0$. In
this case we refer to $U$ and $W$ as {\em orthogonal summands} of
$V$. A bilinear space is {\em indecomposable} if it lacks proper
non-zero orthogonal summands. If $\langle U,U\rangle=0$ then $U$
is {\em totally isotropic}.

For a subspace $U$ of $V$, let $$L(U)=\{v\in V\,\vert\,\langle
v,U\rangle=0\},\quad R(U)=\{v\in V\,\vert\,\langle
U,v\rangle=0\}.$$ Here $L(V)$ and $R(V)$ are the {\em left} and
{\em right radicals} of $V$, and $\R(V)=L(V)\cap R(V)$ is the {\em
radical} of $V$. We have $\mathrm{dim}\,L(V)=\mathrm{dim}\,R(V)$,
and we say that $V$ is {\em non-degenerate} whenever this number
is 0. Otherwise $V$ is {\em degenerate}. A degenerate space is
{\em totally degenerate} if all its non-zero orthogonal summands
are degenerate.

We view $L$ and $R$ as operators which assign to each subspace of
$V$ its left and right orthogonal complements, respectively. If
required we will write $L_V$ and $R_V$ for them. We may compound
these operators, denoting by $L^i$ and $R^i$ their respective
$i$-th iterates. By convention, $L^0$ and $R^0$ are the identity
operators. By definition
$$
L(V)\subseteq L^3(V) \subseteq L^5(V) \subseteq\cdots\subseteq
L^4(V)\subseteq L^2(V) \subseteq L^0(V)=V,
$$
and similarly for $R$. We denote by $L_\infty(V)$, $R_\infty(V)$,
$L^\infty(V)$ and $R^\infty(V)$ the subspaces of $V$ at which the
sequences $(L^{2k+1}(V))_{k\geq 0}$, $(R^{2k+1}(V))_{k\geq 0}$,
$(L^{2k}(V))_{k\geq 0}$ and $(R^{2k}(V))_{k\geq 0}$ stabilize,
respectively. We set $V_\infty=L_\infty(V)+R_\infty(V)$ and
${}^\infty V=L^\infty(V)+R^\infty(V)$. By construction both
$V_\infty$ and ${}^\infty V$ are $G$-invariant.

For $r\geq 1$ and $\lambda\in F$, denote by $J_r(\lambda)$ the
lower Jordan block of size $r$ corresponding to the eigenvalue
$\lambda$. Thus $$ J_1(\lambda)=(\lambda),\quad
J_2(\lambda)=\left(\begin{matrix} \lambda & 0\\ 1
& \lambda\end{matrix}\right),\quad J_3(\lambda)=\left(\begin{matrix} \lambda & 0 & 0\\
1 & \lambda & 0\\ 0 & 1 & \lambda
\end{matrix}\right),...
$$
Write $N_r$ for a bilinear space whose underlying form has matrix
$J_r(0)$ relative to some basis. We shall refer to the bilinear
space $N_r$ as a {\em Gabriel block} and to $r$ as its size. We
refer the reader to \cite{DS1,WW} for the following formulation of
a theorem due to P. Gabriel \cite{PG}.

\begin{theorem}\label {gabriel} Let $(V,\vf)$ be a bilinear space over $F$.
Then

(a) $V=V_{\tdeg}\perp V_{\ndeg}$, where $V_{\tdeg}$ is the
orthogonal direct sum of Gabriel blocks and $V_{\ndeg}$ is
non-degenerate (either of them possibly 0).

(b) The sizes and multiplicities of the Gabriel blocks appearing
in $V_{\tdeg}$ are uniquely determined by $V$.

(c) The equivalence class of $V_{\ndeg}$ is uniquely determined by
$V$.

(d) Up to equivalence, the only indecomposable and degenerate
bilinear space of dimension $r\geq 1$ is $N_r$.
\end{theorem}
We refer to $V_{\tdeg}$ and $V_{\ndeg}$ as the {\em totally
degenerate} and {\em non-degenerate parts} of $V$, respectively.
We may write $V_{\tdeg}=V_{\even}\perp V_{\odd}$, where
$V_{\even}$ resp. $V_{\odd}$ is the orthogonal direct sum of
Gabriel blocks of even resp. odd size. We refer to them as the
{\em even} and {\em odd parts} of $V$.

We fix a decomposition
\begin{equation}
\label{cero} V=V_{\odd}\perp V_{\even}\perp V_{\ndeg},
\end{equation}
and identify $G(V_{\odd})$, $G(V_\even)$ and $G(V_{\ndeg})$ with
their image in $G(V)$, obtained by extending via the identity on
the complements exhibited in (\ref{cero}).

While none of $V_{\even}$, $V_{\odd}$, $V_{\ndeg}$ are in general
$G$-invariant (see \cite{DS1}) we know from \cite{DS1} that,
whatever the choices for these are, $V_\infty$ is the only totally
isotropic subspace of $V_{\odd}$ of maximum dimension and
\begin{equation}
\label{inf} {}^\infty V=V_\infty\perp V_{\even}\perp V_{\ndeg}.
\end{equation}

\begin{notn} If $G$ acts on a set $X$ and $Y\subseteq X$ then $G[Y]$ and
$G\{Y\}$ denote the pointwise and global stabilizers of $Y$ in
$G$, respectively.
\end{notn}

\begin{notn} If $Y$ is a subset of $G$ then $<Y>$ denotes the
subgroup of $G$ generated by $Y$.
\end{notn}

\begin{notn} If $W$ is an $F$-vector space and $f_1,...,f_m$ are
vectors in $W$ then their span will be denoted by $(f_1,...,f_m)$.
\end{notn}

\begin{notn}  The
{\em transpose} of $\phi\in \Bil(V)$ is the bilinear form
$\phi'\in\Bil(V)$, defined by
$$
\phi'(v,w)=\phi(w,v), \quad v,w\in V.
$$
The transpose of a matrix $A$ will be denoted by $A'$.
\end{notn}

\section{Lemmata}

We fix a decomposition
\begin{equation}
\label{uno} V_{\odd}=V_1\perp V_2\perp\cdots\perp V_t,
\end{equation}
where each $V_i$ is the orthogonal direct sum of $m_i$ bilinear
subspaces, each of which is isomorphic to a Gabriel block of size
$2s_i+1$, with $s_1>s_2>\cdots>s_t$. By means of the decomposition
(\ref{uno}) we may identify each $G(V_i)$ with its image in
$G(V_{\odd})$.

We have
\begin{equation}
\label{vik} V_i=V^{i,1}\perp V^{i,2}\perp\cdots \perp V^{i,m_i},
\end{equation}
where each $V^{i,p}$, $1\leq p\leq m_i$, has a basis
$$
e^{i,p}_1,...,e^{i,p}_{2s_i+1}
$$
relative to which the matrix of $\vf$ is equal to $J_{2s_i+1}(0)$.
We shall consider the basis $\B$ of $V_\odd$, defined by
\begin{equation}
\label{bi} \B=\{e^{i,p}_{k}\,|\,1\leq i\leq t,\; 1\leq p\leq
m_i,\; 1\leq k\leq 2s_i+1\}.
\end{equation}
For $1\leq i\leq t$ let $V_i^\dag$ be the span of $e^{i,p}_{2k}$,
$1\leq p\leq m_i$ and $1\leq k\leq s_i$, and let
$$V_{\odd}^\dag=\underset{1\leq i\leq t}\oplus
V_i^\dag.
$$
Note that $V_{\odd}^\dag$ is a totally isotropic subspace of
$V_{\odd}$ satisfying
\begin{equation}
\label{mas} V_{\odd}=V_\infty\oplus V_{\odd}^\dag.
\end{equation}
There is no loss of generality in considering this particular
subspace, as shown in Lemma \ref{regular} below.

\begin{lem}
\label{incl} $G[V_\infty]\subseteq G[V/{}^\infty V]$.
\end{lem}

\noindent{\it Proof.} Let $g\in G[V_\infty]$, $x\in V$ and $y\in
V_\infty$. Then
$$
\la x-gx,y\ra = \la x,y\ra-\la gx,y\ra=\la x,y\ra-\la
x,g^{-1}y\ra=\la x,y\ra-\la x,y\ra=0.
$$
Since $L(V_\infty)={}^\infty V$, the result follows.

\begin{lem}
\label{step1} $G[{}^\infty V]=G(V_\odd) \cap G[V_\infty] \subseteq
G[V/V_\infty]$.
\end{lem}

\noindent{\it Proof.} Since
$$
L(V_\even\oplus V_\ndeg) \cap R(V_\even\oplus V_\ndeg)=V_\odd,
$$
we have $G[{}^\infty V]\subseteq G(V_\odd)$. For $g\in G[{}^\infty
V]$ and $v\in V_\odd$, by Lemma 3.1 we have
$$gv-v\in {}^\infty V \cap V_\odd = V_\infty.
$$
Hence $G[{}^\infty V]\subseteq G[V/V_\infty]$.

\begin{lem}
\label{regular} The permutation action of $G(V_\odd)$ on the set
of totally isotropic subspaces $W$ of $V_{\odd}$ satisfying
$V_{\odd}=V_\infty\oplus W$ is transitive. In fact, restriction to
$G[\VI]$ yields a regular action.
\end{lem}

\noindent{\it Proof.} Let $W$ and $W'$ be totally isotropic
subspaces of $V_{\odd}$ complementing $V_\infty$. By Lemma
\ref{incl} we have
$$G[{}^\infty V]\cap G\{W\}=G[{}^\infty V]\cap G[V/{}^\infty V]\cap G\{W\}=\,<1>.$$
We next show the existence of $g\in G[\VI]$ satisfying $g(W)=W'$.
The decomposition $V_\odd=V_\infty\oplus W'$ gives rise to a
unique projection $p\in \mathrm{End}_F(V_\odd)$ with image $W'$
and kernel $V_\infty$. Define $g\in\GL(V_\odd)$ by
$$g(v+w)=v+p(w),\quad v\in
V_\infty,w\in W.
$$
Let $u,v\in V_\infty$ and $w,z\in W$. Since $V_\infty$, $W$ and
$W'$ are totally isotropic, and $(p-1)V_\odd\subseteq V_\infty$,
we have
$$
\begin{aligned}
\la g(u+w),g(v+z)\ra &= \la u+pw,v+pz\ra=\la u,pz\ra+\la pw,v\ra\\
&=\la u,(p-1)z+z\ra+\la (p-1)w+w,v\ra=\la u,z\ra+\la w,v\ra\\
&=\la u+w,z\ra+\la u+w,v\ra=\la u+w,v+z\ra.
\end{aligned}
$$
Then $g\in G(V_\odd)$ fixes $V_\infty$ elementwise and sends $W$
to $W'$, which completes the proof.

\begin{lem}
\label{block} Let $W=N_{2s+1}$ be a Gabriel block of odd size
$2s+1$. Let $f_1,...,f_{2s+1}$ be a basis of $W$ relative to which
the underlying bilinear form has basis $J_{2s+1}(0)$.

(a) If $0\leq k\leq s$ then
$$
L^{2k+1}(W)=(f_1,f_3,f_5,...,f_{2k+1})\text{ and }
R^{2k+1}(W)=(f_{2s+1},f_{2s-1},f_{2s-3},...,f_{2(s-k)+1}).
$$

(b) If $k\geq s$ then
$$L^{2k+1}(W)=R^{2k+1}(W)=W_\infty=(f_1,f_3,f_5,...,f_{2s+1}).$$
\end{lem}

\noindent{\it Proof.} This follows easily from the definition of
the operators $L$ and $R$.

\begin{lem}
\label{block2} Let $W=N_{2s}$ be a Gabriel block of even size
$2s$. Let $f_1,...,f_{2s}$ be a basis of $W$ relative to which the
underlying bilinear form has basis $J_{2s}(0)$.

(a) If $0\leq k\leq s-1$ then
$$
L^{2k+1}(W)=(f_1,f_3,f_5,...,f_{2k+1})\text{ and }
R^{2k+1}(W)=(f_{2s},f_{2s-2},f_{2s-4},...,f_{2(s-k)}).
$$

(b) If $k\geq s-1$ then
$$L^{2k+1}(W)=(f_1,f_3,f_5,...,f_{2s-1})=L_\infty(W),\;
R^{2k+1}(W)=(f_{2s},f_{2s-2},f_{2s-4},...,f_2)=R_\infty(W).$$

(c) $W=L_\infty(W)\oplus R_\infty(W)$.
\end{lem}

\noindent{\it Proof.} (a) and (b) follow easily from the
definition of the operators $L$ and $R$, and (c) is consequence of
(b).

\begin{lem}
\label{LandR} If $V=U\perp W$ then
$$
L^k_V(V)=L^k_U(U)\perp L^k_W(W)\text{ and }R^k_V(V)=R^k_U(U)\perp
R^k_W(W),\quad k\geq 1.
$$
\end{lem}

\noindent{\it Proof.} This follows easily from the definition of
the operators $L$ and $R$.

\begin{lem}
\label{tower} Let $1\leq i \leq t$ and  $0\leq k,l$.

(a) If  $k,l\leq s_i$. Then a basis for $L^{2k+1}(V)\cap
R^{2l+1}(V)\cap V_i$ is formed by all $e_{2c+1}^{i,p}$, if any,
such that $1\leq p\leq m_i$  and $s_i-l\leq c\leq k$.

(b) If $k>s_i$ (resp. $l>s_i$) then a basis for $L^{2k+1}(V)\cap
R^{2l+1}(V)\cap V_i$ is formed by all $e_{2c+1}^{i,p}$ such that
$1\leq p\leq m_i$, $0\leq c\leq s_i$, and $s_i-l\leq c$ (resp.
$c\leq k$).
\end{lem}

\noindent{\it Proof.} This follows from Lemmas \ref{block} and
\ref{LandR} by means of the decompositions (\ref{cero}),
(\ref{uno}) and (\ref{vik}).

\begin{lem}
\label{toohigh} Let $k,l\geq 0$ and $1\leq i\leq t$. Then
$$
L^{2k+1}(V)\cap R^{2l+1}(V)\cap V_i\neq (0)
$$
if and only if $k+l\geq s_{i}$.
\end{lem}

\noindent{\it Proof.} This follows from Lemma \ref{tower}.

\begin{lem}
\label{tres} Let $1\leq i\leq t$ and $0\leq j,k$. Suppose $i+j\leq
t$ and $k\leq s_i$. Then
$$
L(V)\cap R^{2(s_i-k)+1}(V)\cap V_{i+j}\neq (0)
$$
if and only if $s_i-k\geq s_{i+j}$.
\end{lem}

\noindent{\it Proof.} This is a particular case of Lemma
\ref{toohigh}.

\begin{lem}
\label{cuatro} Let $1\leq i\leq t$ and $0\leq k\leq s_i$. Then
$$
L^{2k+1}(V)\cap R^{2(s_i-k)+1}(V)\cap
V_{i}=(e_{2k+1}^{i,1},...,e_{2k+1}^{i,m_i}).
$$
\end{lem}

\noindent{\it Proof.} This is a particular case of Lemma
\ref{tower}.

\begin{defn} Consider the subspaces of $V_\infty$ defined as follows:
$$
V(i)=\underset{i\leq j\leq t}{\oplus} (V_j)_\infty,\quad 1\leq
i\leq t
$$
and set $V(i)=0$ for $i>t$.
\end{defn}

\begin{lem}
\label{cola} If $1\leq i\leq t$ then
$$
V(i)=\underset{0\leq k\leq s_i }\sum L^{2k+1}(V)\cap
R^{2(s_i-k)+1}(V).
$$
\end{lem}

\noindent{\it Proof.} By virtue of Lemmas \ref{block2},
\ref{LandR} and \ref{toohigh}, and the decompositions
(\ref{cero}), (\ref{uno}) and (\ref{vik}) it follows that the
right hand side is contained in $V(i)$. By Lemma \ref{cuatro}, if
$i\leq j\leq t$ then
$$
(V_j)_\infty=\underset{0\leq k\leq s_j }\sum L^{2k+1}(V)\cap
R^{2(s_j-k)+1}(V)\cap V_j\subseteq \underset{0\leq k\leq s_i }\sum
L^{2k+1}(V)\cap R^{2(s_i-k)+1}(V),$$ as required.

\begin{lem}
\label{inv} The subspaces $V(i)$ are $G$-invariant.
\end{lem}

\noindent{\it Proof.} This follows from Lemma \ref{cola}.

\begin{lem}
\label{radi} Let $W=N_{2s+1}$ be a Gabriel block of odd size
$2s+1$. Let $f_1,...,f_{2s+1}$ be a basis of $W$ relative to which
the underlying bilinear form, say $\phi$, has basis $J_{2s+1}(0)$.
Then $\R(\phi-\phi')=(f_1+f_3+\cdots+f_{2s-1}+f_{2s+1})$.
\end{lem}

\noindent{\it Proof.} Clearly the vector
$f_1+f_3+\cdots+f_{2s-1}+f_{2s+1}$ belongs to the radical of
$\phi-\phi'$. Since the nullity of the matrix
$J_{2s+1}(0)-J_{2s+1}(0)'$ is equal to one, the result follows.

\begin{notn} For each $1\leq i\leq t$ and each $1\leq p\leq m_i$
let
$$
E^{i,p}=e^{i,p}_{1}+e^{i,p}_{3}+\cdots+e^{i,p}_{2s_i+1}.
$$
\end{notn}

\begin{lem}
\label{rad} If $1\leq i\leq t$ then
$$
\R(\vf-\vf')\cap V_i=(E^{i,1},...,E^{i,m_i}).
$$
\end{lem}

\noindent{\it Proof.} By Lemma \ref{radi} we have
$$
\R(\vf-\vf')\cap V_i=\underset{1\leq p\leq m_i}\oplus
\R(\vf-\vf')\cap V^{i,p}=\underset{1\leq p\leq m_i}\oplus
(E^{i,p})=(E^{i,1},...,E^{i,m_i}).
$$

\begin{defn}
\label{bert1} Let $1\leq i$; $0\leq k,j$; $1\leq p,q$. Suppose
$i+j\leq t$; $\,k\leq s_i-s_{i+j}$; $1\leq p\leq m_i$; $1\leq
q\leq m_{i+j}$; $p\neq q$ if $j=0$. Consider the 1-parameter
subgroup of $G(V^{i,p}\perp V^{i+j,q})$ -or simply $G(V^{i,p})$ if
$j=0$- formed by all $g^{i,i+j,p,q}_{2k+1,y}\in G(V_{\odd})$, as
$y$ runs through $F$, defined as follows.

For ease of notation we replace $g^{i,i+j,p,q}_{2k+1,y}$ by $g$;
$s_i$ by $s$; $e^{i,p}_1,...,e^{i,p}_{2s+1}$ by
$e_1,...,e_{2s+1}$; $s_{i+j}$ by $d$; and
$e^{i+j,q}_1,...,e^{i+j,q}_{2d+1}$ by $f_1,...,f_{2d+1}$. If $v\in
V_\odd$ then
$g$ fixes all basis vectors of (\ref{bi}) not listed below and
$$
g(e_{2k+1})=e_{2k+1}+yf_1,\quad g(f_2)=f_2-ye_{2k+2},
$$
$$
g(e_{2k+3})=e_{2k+3}+yf_3,\quad g(f_4)=f_4-ye_{2k+4},
$$
$$
\vdots
$$
$$
g(e_{2(k+d)-1})=e_{2(k+d)-1}+yf_{2d-1},\quad
g(f_{2d})=f_{2d}-ye_{2(k+d)},
$$
$$
g(e_{2(k+d)+1})=e_{2(k+d)+1}+yf_{2d+1}.
$$
To see that $g$ indeed belongs to $G(V_\odd)$ it suffices to
verify that the matrices of $\vf$ relative to the bases $\B$ and
$\vf(\B)$ are equal. This is a simple computation involving basis
vectors from at most two Gabriel blocks, and we omit it.
\end{defn}

\begin{defn}
\label{bert2}
 For $1\leq i\leq t$, $1\leq p\leq m_i$ consider the
1-parameter subgroup of $G(V^p_i)$ formed by all $g^{i,p}_{x}\in
G(V_{\odd})$, as $x$ runs through $F^*$, defined as follows.

For ease of notation we replace $g^{i,p}_x$ by $g$; $s_i$ by $s$;
and $e^{i,p}_1,e^{i,p}_2,...,e^{i,p}_{2s+1}$ by
$e_1,e_2,...,e_{2s+1}$. If $v\in V^p_i$ then
$g$ fixes all basis vectors of (\ref{bi}) not listed below and
$$
g(e_1)=xe_1, g(e_2)=x^{-1}e_2,...,$$ $$g(e_{2s-1})=xe_{2s-1},
g(e_{2s})=x^{-1}e_{2s}, g(e_{2s+1})=xe_{2s+1}.
$$
In this this case one easily verifies that $g\in G(V_\odd)$.
\end{defn}

\begin{lem}
\label{basis} Suppose $W$ is an $F$-vector space with a basis
$$
f_1^1,...,f_1^m,f_2^1,...,f_2^m,...,f_s^1,...,f_s^m.
$$
For each $1\leq p\leq m$ let
$$
E^p=f_1^p+f_2^p+\cdots+f_s^p.
$$
Suppose $g\in\mathrm{End}_F(W)$ preserves each of the
$m$-dimensional subspaces $(f_k^1,...,f_k^m)$, where $1\leq k\leq
s$, and also the $m$-dimensional subspace $(E^1,...,E^m)$. Suppose
further that $g$ fixes all vectors $f_1^1,...,f_1^m$. Then $g=1$.
\end{lem}

\noindent{\it Proof.} From the invariance of the subspaces
$(f_k^1,...,f_k^m)$ we have
$$
g(f_k^p)=\underset{1\leq q\leq m}\sum a_k^{p,q}f_k^q,
$$
where $a_k^{p,q}\in F$. Since $g$ fixes $f_1^1,...,f_1^m$
$$
a_1^{p,q}=\delta_{p,q}.
$$
As $g$ is linear
$$
g(E^p)=\underset{1\leq k\leq s}\sum\; \underset{1\leq q\leq m}\sum
a_k^{p,q}f_k^q=\underset{1\leq q\leq m}\sum \;\underset{1\leq
k\leq s}\sum a_k^{p,q} f_k^q.
$$
But by the invariance of $(E^1,...,E^m)$ we also have
$$
g(E^p)=\underset{1\leq q\leq m}\sum b^{p,q}E^q=\underset{1\leq
q\leq m}\sum \;\underset{1\leq k\leq s}\sum b^{p,q} f_k^q,
$$
where $b^{p,q}\in F$. Therefore $a_k^{p,q}$ is independent of $k$,
and in particular $$a_k^{p,q}=a_1^{p,q}=\delta_{p,q},$$ as
required.

\section{The split extension $1\to N\to G\to G/N\to 1$}
\label{lo}

\begin{defn} For $j\geq 1$ consider the subgroup
$N_j$ of $G$ defined by
$$N_j=\underset{1\leq i\leq
t}\cap G[V(i)/V(i+j)],
$$
and set $N=N_1$. Each $N_j$ is normal due to Lemma \ref{inv}.
\end{defn}
Note that
$$N=N_1\supseteq N_2\supseteq\cdots \supseteq
N_t=G[V_\infty],\quad N_j=G[V_\infty],\text{ if } j\geq t.$$

\begin{defn} For $1\leq i\leq t$ let
$E_i$ be the subgroup of $G(V_i)$ generated by all
$g^{i,i,p,q}_{1,y}$ and all $g^{i,p}_x$. Let $E$ be the subgroup
of $G(V_{\odd})$ generated by all $E_i$, $1\leq i\leq t$. Let
$$
E_i'=G(V_i)\cap G\{V_i^\dag \},\quad 1\leq i\leq t,
$$
and consider the internal direct product
$$E'=\underset{1\leq i\leq t}\Pi E_i'.
$$
\end{defn}

\begin{defn}
\label{desi}
 For $1\leq i\leq t$ and $0\leq k\leq s_i$ consider
the $FG$-submodule $S_{2k+1}^i$ of $V(i)/V(i+1)$, defined by
$$
S_{2k+1}^i=\left(L^{2k+1}(V)\cap R^{2(s_i-k)+1}(V)\cap
V(i)+V(i+1)\right)/V(i+1).
$$
\end{defn}
We know from Lemma \ref{cuatro} that
$$
S_{2k+1}^i=\left((e_{2k+1}^{i,1},...,e_{2k+1}^{i,m_i})\oplus
V(i+1)\right)/V(i+1).
$$
This yields the following decomposition of $FG$-modules
$$
V(i)/V(i+1)=\underset{0\leq k\leq s_i}\oplus S_{2k+1}^i.
$$

\begin{thm}
\label{E} The canonical map
\begin{equation}
\label{epi} G\to \underset{1\leq i\leq t}\Pi \GL(S^i_1)\cong
\underset{1\leq i\leq t}\Pi \GL_{m_i}(F)
\end{equation}
is a split group epimorphism with kernel $N$. Moreover,
\begin{equation}
\label{LA} E_i'=E_i\cong\GL_{m_i}(F)\text{ for all }1\leq i\leq t,
\end{equation}
$$
E'=E,
$$
and
$$G=N\rtimes E.
$$
\end{thm}

\noindent{\it Proof.} The above description of $S^i_1$ combined
with Definitions \ref{bert1} and \ref{bert2} yield that each map
$$
E_i\to \GL(S^i_1),\quad 1\leq i\leq t.
$$
is surjective, whence the map (\ref{epi}) is also surjective.

To see that the kernel of (\ref{epi}) is $N$ we apply Lemmas
\ref{cuatro}, \ref{rad} and \ref{basis}. Indeed, if $g\in G$ is in
the kernel of (\ref{epi}) then the $G$-invariance of
$\R(\vf-\vf')\cap V(i)+V(i+1)/V(i+1)$ and all $S_{2k+1}^i$, $0\leq
k\leq s_i$, along with the fact that $g$ acts trivially on
$S_1^i$, imply that $g$ acts trivially on $V(i)/V(i+1)$ for all
$1\leq i\leq t$, as required.

It follows that $G=NE$. But by definition $E\subseteq E'$ and
$E'\cap N=1$. Therefore $E'=E$, $G=N\rtimes E$, and $E_i\cong
\GL_{m_i}(F)$ for all $1\leq i\leq t$. Since it is obvious that
$E$ is the internal direct product of the $E_i\subseteq E_i'$, the
proof is complete.

\section{Irreducible constituents of $V$ as an $FG$-module}
\label{comp}

The series
$$
0\subseteq V_\infty\subseteq {}^\infty V\subseteq V
$$
reduces the search of irreducible constituents of the $FG$-module
$V$ to that of the factors
$$
V_\infty,\quad {}^\infty V/V_\infty,\quad V/{}^\infty V.
$$
First we consider the factor $V_\infty$.

\begin{thm}
\label{irrvi} Each factor $V(i)/V(i+1)$, $1\leq i\leq t$, of the
series of $FG$-modules
$$
V_\infty=V(1)\supset V(2)\supset\cdots\supset V(t)\supset V(t+1)=0
$$
is equal to the direct sum of $s_i+1$ isomorphic irreducible
$FG$-modules of dimension $m_i$
$$
V(i)/V(i+1)=\underset{0\leq k\leq s_i}\oplus S_{2k+1}^i.
$$
Moreover,
$$
G[S_{2k+1}^i]=N\rtimes \underset{l\neq i}\Pi E_i,
$$
and as a module for
$$
G/G[S^i_{2k+1}]\cong E_i\cong \GL_{m_i}(F),
$$
$S_{2k+1}^i$ is isomorphic to the natural $m_i$-dimensional module
over $F$, namely $F^{m_i}$.
\end{thm}

\noindent{\it Proof.} This is clear from section \S\ref{lo}.

Next we make preliminary remarks about the factor ${}^\infty
V/V_\infty$.

\begin{defn}
\label{eo}
 Let $V^\infty=L^\infty(V)\cap R^\infty(V)$ and
${}_\infty V=L_\infty(V)+R_\infty(V)$.
\end{defn}

By construction these are $FG$-submodules of $V$, and we know from
\cite{DS1} that
$$
{}_\infty V=V_\even\perp V_\infty\text{ and }
V^\infty=V_\ndeg\perp V_\infty.
$$
Observe that ${}^\infty V/V_\infty$ is a bilinear space, with even
and non-degenerate parts equal to
$${}_\infty V/V_\infty\cong V_\even\text{ and }V^\infty/V_\infty\cong V_\ndeg.$$
Since the odd part of ${}^\infty V/V_\infty$ is equal to zero, we
know from \cite{DS1} that even and non-degenerate parts of
${}^\infty V/V_\infty$ are unique, so
$$
G({}^\infty V/V_\infty)=G({}_\infty V/V_\infty)\times G(V^\infty
/V_\infty).
$$
It follows that the canonical map
$$
G(V)\to G({}^\infty V/V_\infty)
$$
is a group epimorphism whose restriction to $G(V_\even)\times
G(V_\ndeg)$ is an isomorphism. Since the kernel of this map is
$G[{}^\infty V/V_\infty]$, whose intersection with
$G(V_\even)\times G(V_\ndeg)$ is trivial, we obtain the
decomposition
$$
G(V)=G[{}^\infty V/V_\infty]\rtimes G(V_\even)\times G(V_\ndeg).
$$
It follows from the above considerations that the study of the
$FG$-module ${}^\infty V/V_\infty$ reduces to the study of the
$FG(V_\even)$-module $V_\even$ and the $FG(V_\ndeg)$-module
$V_\ndeg$.

\section{The split extension $1\!\to\! G[V_\infty]\!\to\! N\!\to\!
N/G[V_\infty]\!\to\! 1$} \label{nginfty}

The very definition of the groups $N_j$ gives
$$
[N_i,N_j]\subseteq N_{i+j},
$$
so $(N_j)_{1\leq j\leq t}$ yields a $G$-invariant descending
central series for $N/G[V_\infty]$. By abuse of language we shall
sometimes say that $(N_j)_{1\leq j\leq t}$ and like series are
central series for $N/G[V_\infty]$.

\begin{notn}
If $g_1,g_2\in G$ then $[g_2,g_1]=g_2^{-1}g_1^{-1}g_2g_1$. If
$n>2$ and $g_1,...,g_{n-1},g_n\in G$ then
$[g_n,g_{n-1},...,g_1]=[g_n,[g_{n-1},...,g_1]]$.
\end{notn}

\begin{thm}
\label{nilpdegree} The nilpotency class of $N/G[V_\infty]$ is
$t-1$.
\end{thm}

\noindent{\it Proof.} By the above comments the nilpotency class
of $N/G[V_\infty]$ is at most $t-1$. If $t>1$ then
$$
[g^{t-1,t,1,1}_{1,1},...,g^{2,3,1,1}_{1,1},g^{1,2,1,1}_{1,1}]\neq
1,
$$
so the result follows.

The series $(N_j)_{1\leq j\leq t}$ needs to be refined in order to
obtain sharper results on the structure of $N/G[V_\infty]$.

\subsection{Generators for nilpotent group $N/G[V_\infty]$}

\begin{defn} For $k\geq 0$ define the normal subgroup
$M_{2k+1}$ of $G$ by
$$
M_{2k+1}=G[L^{2k+1}(V)\cap V_\infty]\cap N.
$$
We further define $M_{-1}=N$.
\end{defn}
Note that
$$
N=M_{-1}\supseteq M_1\supseteq M_3\supseteq \cdots\supseteq
M_{2s_1+1}=G[V_\infty],
$$
with
$$
M_{2k+1}=G[V_\infty],\quad k\geq s_1.
$$
We use the $G$-invariant series $(M_{2k-1})_{0\leq k}$ to refine
the $G$-invariant decreasing central series $(N_j)_{1\leq j}$ for
$N/G[V_\infty]$, obtaining the $G$-invariant decreasing central
series for $N/G[V_\infty]$
\begin{equation}
\label{nbv} N_{j,2k-1}=(N_j\cap M_{2k-1})N_{j+1},\quad 0\leq k,
1\leq j.
\end{equation}
We have $N_{j,2s_1+1}=N_{j+1}=N_{j+1,-1}$ and
$$
N_1=N_{1,-1}\supseteq N_{1,1}\supseteq
N_{1,3}\supseteq\cdots\supseteq N_{1,2s_1-1}\supseteq
N_2\supseteq\cdots
$$
$$
N_{t-1}=N_{t-1,-1}\supseteq N_{t-1,1}\supseteq\cdots \supseteq
N_{t-1,2s_1-1}\supseteq N_{t-1,2s_1+1}=N_t=1.
$$

\begin{thm}
\label{sent-to-L}
 Let $k\geq 0$. Then
$$
M_{2k-1}\subseteq G[L^{2k+1}(V)\cap V_\infty/L(V)\cap V_\infty].
$$
\end{thm}

\noindent{\it Proof.} We may assume $k\geq 1$, for otherwise the
result is trivial. Since
\begin{equation}
\label{dec} L^{2k+1}(V)\cap V_\infty=L^{2k+1}(V)\cap
V_1\oplus\cdots\oplus L^{2k+1}(V)\cap V_t,
\end{equation}
it suffices to show
\begin{equation}
\label{vi} (g-1)L^{2k+1}(V)\cap V_i\subseteq L(V)\cap
V_\infty,\quad g\in M_{2k-1}, 1\leq i\leq t.
\end{equation}
Fix $i$, $1\leq i\leq t$. If $k>s_i$ then $L^{2k+1}(V)\cap V_i=
L^{2k-1}(V)\cap V_i$, so (\ref{vi}) holds. Suppose $k\leq s_i$.
Then
\begin{equation}
\label{dec2} L^{2k+1}(V)\cap V_i=L^{2k-1}(V)\cap V_i\oplus
L^{2k+1}(V)\cap R^{2(s_i-k)+1}(V)\cap V_i,
\end{equation}
so (\ref{vi}) is equivalent to
\begin{equation}
\label{vi2} (g-1)L^{2k+1}(V)\cap R^{2(s_i-k)+1}(V)\cap V_i
\subseteq L(V)\cap V_\infty,\quad g\in M_{2k-1}.
\end{equation}
By Lemma \ref{cuatro} a basis for $L^{2k+1}(V)\cap
R^{2(s_i-k)+1}(V)\cap V_i$ is given by $(e_{2k+1}^{i,p})_{1\leq
p\leq m_i}$. Let $g\in M_{2k-1}$ and fix $p$, $1\leq p\leq m_i$.
We are reduced to show that $g$ fixes $e_{2k+1}^{i,p}$ modulo
$L(V)\cap V_\infty$. Since $g\in N$, we have
$$
g(e^{i,p}_{2k+1})=e^{i,p}_{2k+1}+z,
$$
where $z\in L^{2k+1}(V)\cap V(i+1)$. Suppose $z\notin L(V)\cap
V(i+1)$. Then
$$
\langle ge^{i,p}_{2k+1}, e^{l,q}_{2c}\rangle\neq 0$$ for some
$t\geq l\geq i+1$, $1\leq q\leq m_l$ and $1\leq c\leq s_l$. Then
$$\langle e^{i,p}_{2k+1}, g^{-1}e^{l,q}_{2c}\rangle\neq 0,
$$
so $g^{-1}e^{l,q}_{2c}$ has non-zero coefficient in
$e^{i,p}_{2k}$. But $g\in M_{2k-1}$ and $l>i$, so
$$
0\neq \langle g^{-1}e^{l,q}_{2c}, e^{i,p}_{2k-1}\rangle= \langle
e^{l,q}_{2c}, ge^{i,p}_{2k-1}\rangle= \langle e^{l,q}_{2c},
e^{i,q}_{2k-1}\rangle=0,
$$
a contradiction.

\begin{defn} Let $1\leq j$ and $0\leq k$. Set
$$
I(j,k)=\{i\geq 1\,|\, 1\leq i\leq t-j\text{ and }k\leq
s_i-s_{i+j}\}.
$$
Note that $I(j,k)=\emptyset$ if $j\geq t$ or $k>s_1$. For $i\geq
1$ we set $X(i,j,2k-1)=\emptyset$ of $i\notin I(j,k)$ and
otherwise
$$
X(i,j,2k-1)= \{g^{i,i+j,p,q}_{2k+1,y}\,|\, 1\leq p\leq m_i, 1\leq
q\leq m_{i+j}, y\in F\}.
$$
We further define
$$
X(j,2k-1)=\underset{i\geq 1}\cup X(i,j,2k-1),\;
X(j)=\underset{k\geq 0}\cup X(j,2k-1),\; X=\underset{j\geq 1}\cup
X(j).
$$
\end{defn}

\begin{thm}
\label{important} Let $k\geq 0$ and $j\geq 1$. Then
$X(j,2k-1)\subseteq M_{2k-1}\cap N_j$ and the quotient group
$$(M_{2k-1}\cap N_j)M_{2k+1}/(M_{2k-1}\cap N_{j+1})M_{2k+1}$$
is generated by the classes of all elements in $X(j,2k-1)$. In
particular, this quotient is trivial if $I(j,k)=\emptyset$ (the
converse is true and proved in Theorem \ref{G-action} below).
\end{thm}

\noindent{\it Proof.} Clearly $X(j,2k-1)\subseteq M_{2k-1}\cap
N_j$. Let $g\in M_{2k-1}\cap N_j$. We claim that for each $i$,
$1\leq i\leq t$, there exists
\begin{equation}\
\label{hi} h_i\in < X(i,j,2k-1)> < \underset{r>j}\cup X(i,r,2k-1)>
\end{equation}
such that $h_ig$ is the identity on $L^{2k+1}(V)\cap V_i$.

To prove our claim we fix $i$, $1\leq i\leq t$. If $k>s_i$ we take
$h_i=1$. Suppose $k\leq s_i$. In view of $g\in M_{2k-1}$ and the
decomposition (\ref{dec2}), it suffices to choose $h_i$ as in
(\ref{hi}), so that $h_ig$ fixes every basis vector
$e^{i,p}_{2k+1}$, $1\leq p\leq m_i$, of $L^{2k+1}(V)\cap
R^{2(s_i-k)+1}(V)\cap V_i$. By Theorem \ref{sent-to-L} and the
fact that $g\in N_j$, for each $1\leq p\leq m_i$ we have
$$
g(e^{i,p}_{2k+1})=e^{i,p}_{2k+1}+z,
$$
where $z\in L(V)\cap R^{2(s_i-k)+1}(V)\cap V(i+j)$. If $i+j>t$
then $z=0$ and if $i+j\leq t$ but $k>s_i-s_{i+j}$ then $z=0$ as
well by Lemma \ref{tres}. In both cases we take $h_i=1$.
Otherwise, again by Lemma \ref{tres}, we may write
$$
z=z_{i+j}+\cdots+z_{i+j+l},
$$
where $l\geq 0$, $z_{i+j+b}\in L(V)\cap R^{2(s_i-k)+1}(V)\cap
V_{i+j+b}$ and $k\leq s_i-s_{i+j+b}$ for all $0\leq b\leq l$.

It is know clear from the very definition of the
$g^{i,i+j+b,p,q}_{2k+1,y}$ that we may choose $h_i$ as in
(\ref{hi}) so that $h_ig$ is the identity on each
$e^{i,p}_{2k+1}$.

By construction $h_i$ is the identity on $\underset{l\neq i}\oplus
V_l$. Therefore our claim and the decomposition (\ref{dec}) imply
that $h_1\cdots h_tg$ is the identity on $L^{2k+1}(V)\cap
V_\infty$, i.e. $h_1\cdots h_tg\in M_{2k+1}$. But if $j<r$ then
$X(i,r,2k-1)\subseteq M_{2k-1}\cap N_{j+1}$, so the class of $g$
is equal to the product of the classes of elements from
$X(j,2k-1)$, as required.

\begin{thm}
\label{butterfly} Let $k\geq 0$ and $j\geq 1$. Then
$X(j,2k-1)\subseteq M_{2k-1}\cap N_j$ and the quotient group
$$(M_{2k-1}\cap N_j)N_{j+1}/(M_{2k+1}\cap N_j)N_{j+1}$$
is generated by the classes of all elements in $X(j,2k-1)$. In
particular, this quotient is trivial if $I(j,k)=\emptyset$ (the
converse is true and proved in Theorem \ref{G-action} below).
\end{thm}

\noindent{\it Proof.} This follows from Theorem \ref{important}
and the Butterfly Lemma.

\begin{thm}
\label{generation} Let $1\leq j<j'\leq t$. Then $X(l)\subseteq
N_j$ for all $j\leq l<j'$ and $N_j/N_{j'}$ is generated the
classes of all these elements. In particular $N/G[V_\infty]$ is
generated by the classes of all elements in $X$.
\end{thm}

\noindent{\it Proof.} This follows from Theorem \ref{butterfly}
via the series $N_{j,2k-1}$.

\begin{defn} Let $U$ be the subgroup of $G(V_{\odd})\cap N$ generated by
$X$ and let
$$
U'=N\cap G(V_{\odd})\cap G[V_1^\dag]\cap G[V_1^\dag\oplus
V_2^\dag/V_1^\dag]\cap\cdots\cap G[V_1^\dag\oplus\cdots \oplus
V_t^\dag/V_1^\dag\oplus\cdots\oplus V_{t-1}^\dag].
$$
\end{defn}

It is clear that $N/G[V_\infty]$ is a unipotent subgroup of
$V_\infty$. The next result shows that $N=G[V_\infty]\rtimes U$,
where $U$ is unipotent in $V$.

\begin{thm}
\label{U} $U'=U$ is unipotent and $N=G[V_\infty]\rtimes U$.
\end{thm}

\noindent{\it Proof.} From Theorem \ref{generation} we infer
$N=G[V_\infty]U$. The definition of $X$ yields $U\subseteq U'$.
But by Lemma \ref{incl}
$$
G[V_\infty]\cap U'\subseteq G[V_\infty]\cap G[V/{}^\infty V]\cap
U'=1.
$$
It follows that $U=U'$ and $N=G[V_\infty]\rtimes U$.

Finally, if $g\in U$ then $g-1\vert_{V_{\odd}^\dag}$ is nilpotent
by the very definition of $U'$, and $g-1\vert_{{}^\infty V}$ is
nilpotent since $g\in N$. Thus $g$ is unipotent.

\subsection{Irreducible constituents of the $FG$-module
$V/{}^\infty V$}

We have accumulated enough information to determine the
irreducible constituents of the $FG$-module $V/{}^\infty V$.

\begin{defn} Let $\t=t$ if $s_t>0$ (i.e. $\R(V)= 0$) and
$\t=t-1$ if $s_t=0$ (i.e. $\R(V)\neq 0$).
\end{defn}

\begin{defn} For $1\leq i\leq \t$ let
$$(i)V=V_1^\dag\oplus\cdots\oplus V_i^\dag\oplus \VI$$
and set
$$
(0)V=\VI.
$$
\end{defn}

\begin{lem}
\label{inb} If $0\leq i\leq \t$ then $(i)V$ is an $FG$-submodule
of $V$.
\end{lem}

\noindent{\it Proof.} We may assume that $i\geq 1$. From the
identity $N=G[V_\infty]\rtimes U$, the inclusion
$G[V_\infty]\subseteq G[\WI]$ of Lemma \ref{incl}, and the
characterization of $U$ given in Theorem \ref{U} we infer that
$(i)V$ is preserved by $N$. The very definition of $E$ and the
identity $G=N\rtimes E$ of Theorem \ref{E} allow us to conclude
that $(i)V$ is in fact $G$-invariant.

\begin{lem}
\label{inb2} If $1\leq i\leq \t$ then $(i)V/(i-1)V$ is an
$FG$-module acted upon trivially by $N$.
\end{lem}

\noindent{\it Proof.} The characterization of of $U$ given in
Theorem \ref{U} shows that $U$ acts trivially on $(i)V/(i-1)V$,
while Lemma \ref{incl} shows that $G[V_\infty]$ also acts
trivially on $(i)V/(i-1)V$. Since $N=G[V_\infty]\rtimes U$, the
result follows.


\begin{thm}
\label{irrvI} Each factor $(i)V/V(i-1)$, $1\leq i\leq \t$, of the
series of $FG$-modules
$$
\VI=(0)V\subset (1)V\subset\cdots\subset (\t-1)V\subset (\t)V=V
$$
is isomorphic to the direct sum of $s_i$ isomorphic irreducible
$FG$-modules of dimension $m_i$, namely
$$
Q_{2k}^i=\left((e_{2k}^{i,1},...,e_{2k}^{i,m_i})\oplus
(i-1)V\right)/(i-1)V,$$ where $1\leq k\leq s_i$. Moreover,
$$
G[Q_{2k}^i]=N\rtimes \underset{l\neq i}\Pi E_i,
$$
and as a module for
$$
G/G[Q^i_{2k}]\cong E_i\cong \GL_{m_i}(F),
$$
$Q_{2k}^i$ is isomorphic to the natural $m_i$-dimensional module
over $F$, namely $F^{m_i}$.
\end{thm}

\noindent{\it Proof.} This follows from Lemma \ref{inb2} and
(\ref{LA}).

\subsection{A refined $G$-invariant descending central series for
$N/G[V_\infty]$}

In this section we construct a $G$-invariant descending central
series for $N/G[V_\infty]$ each of whose factors is naturally an
irreducible $FG$-module, whose isomorphism type we explicitly
determine, all of which are connected to the irreducible
constituents of the $FG$-module $V_\infty$, as described in
Theorem \ref{irrvi}

\begin{thm}
\label{vectorspace} For each $j\geq 1$ there is a canonical group
embedding
$$
N_j/N_{j+1}\to \underset{1\leq i\leq t-1}\oplus
\mathrm{Hom}_F(V(i)/V(i+1),V(i+j)/V(i+j+1)),
$$
whose image is an $F$-vector subspace of the codomain. By
transferring this $F$-vector space structure to $N_j/N_{j+1}$ the
above map becomes an embedding of $FG$-modules.
\end{thm}

\noindent{\it Proof.} Recall first of all that since $V(i)/V(i+1)$
and $V(i+j)/V(i+j+1)$ are $FG$-modules, so is
$\mathrm{Hom}_F(V(i)/V(i+1),V(i+j)/V(i+j+1))$ in a natural manner.

Define the map $N_j/N_{j+1}\to
\mathrm{Hom}_F(V(i)/V(i+1),V(i+j)/V(i+j+1))$, $1\leq i\leq t-1$,
by $[g]\mapsto g_i$, where $g\in N_j$ and
$$
g_i(v+V(i+1))=(g-1)(v)+V(i+j+1),\quad v\in V(i).
$$
Then let $N_j/N_{j+1}\to \underset{1\leq i\leq t-1}\oplus
\mathrm{Hom}_F(V(i)/V(i+1),V(i+j)/V(i+j+1))$ be defined by
$$
[g]\mapsto (g_1,...,g_{t-1}),\quad g\in N_j.
$$
The definitions of all objects involved and the identity
$$
gh-1=(g-1)(h-1)+(h-1)+(g-1),\quad g,h\in G
$$
show that our map is a well-defined group monomorphism which is
compatible with the action of $G$ on both sides.

It remains to show that the image of our map is an $F$-subspace of
the codomain. By Theorem \ref{generation} $N_j/N_{j+1}$ is
generated by all $gN_{j+1}$ as $g$ runs through $X(j)$. Since our
map is a group homomorphism, it suffices to show that $k(g-1)+1\in
N_j$ for all $g\in X(j)$ and $k\in F$. But Definition \ref{bert1}
makes this clear, so the proof is complete.

\begin{note} Let $k\geq 0$ and $j\geq 1$. Notice that the group
$N_{j,2k-1}/N_{j,2k+1}$ is a section of the $FG$-module
$N_j/N_{j+1}$ of Theorem \ref{vectorspace}, and as such inherits a
natural structure of $FG$-module.
\end{note}

\begin{defn}
\label{doce}
 Let $k\geq 0$ and $j\geq 1$. For each $1\leq
i\leq t$ define
$$
N_{i,j,2k-1}=\{g\in N_{j,2k-1}\,\vert\, (g-1)L^{2k+1}(V)\cap
V(l)\subseteq V(l+j+1)\text{ for all }1\leq l\leq t, l\neq i\}.
$$
\end{defn}

Note that for each $1\leq i\leq t$, $N_{i,j,2k-1}$ is a normal
subgroup of $G$ containing $N_{j,2k+1}$. As a section of
$N_j/N_{j+1}$, the the group $N_{i,j,2k-1}/N_{j,2k+1}$ is also an
$FG$-module.

Recall at this point the meaning of the $FG$-modules $S_{2k+1}^i$,
as given in Definition \ref{desi}.

\begin{thm}
\label{G-action} Let $k\geq 0$ and $j\geq 1$. For each $i\in
I(j,k)$ there is a canonical isomorphism of $FG$-modules
\begin{equation}
\label{ula} Y_i:N_{i,j,2k-1}/N_{j,2k+1}\to \mathrm{Hom}_F
(S^i_{2k+1},S_1^{i+j}).
\end{equation}
Moreover, we have $X(i,j,2k-1)\subseteq N_{i,j,2k-1}$, and
$N_{i,j,2k-1}/N_{j,2k+1}$ is generated by the classes of all
elements in $X(i,j,2k-1)$.

There is a canonical isomorphism of $FG$-modules
$$
Y:(N_j\cap M_{2k-1})N_{j+1}/(N_j\cap M_{2k+1})N_{j+1}\to
\underset{i\in I(j,k)}\oplus\mathrm{Hom}_F (S^i_{2k+1},S_1^{i+j}).
$$
induced by the $Y_i$. The dimension of both of these modules, say
$d_{j,2k-1}$, is equal to $d_{j,2k-1}=\underset{i\in I(j,k)}\sum
m_im_{i+j}$. Moreover, we have
\begin{equation}
\label{nice} (N_j\cap M_{2k-1})N_{j+1}/(N_j\cap
M_{2k+1})N_{j+1}=\underset{i\in I(j,k)}\Pi
N_{i,j,2k-1}/N_{j,2k+1}.
\end{equation}
\end{thm}

\noindent{\it Proof.} Let $i\in I(j,k)$. In the spirit of Theorem
\ref{vectorspace} we consider the map
$$Y_i: (N_j\cap M_{2k-1})N_{j+1}/(N_j\cap
M_{2k+1})N_{j+1} \to \mathrm{Hom}_F (S^i_{2k+1},S_1^{i+j})$$ given
by $[g]\mapsto g_i$, where $g\in (N_j\cap M_{2k-1})N_{j+1}$ and
$$
g_i(v+V(i+1))=(g-1)(v)+V(i+j+1),
$$
for all $v\in L^{2k+1}(V)\cap R^{2(s_i-k)+1}(V)\cap V(i)+V(i+1)$.

{\sc Step I:} $Y_i$ is a well-defined homomorphism of
$FG$-modules.

Let $g\in (N_j\cap M_{2k-1})N_{j+1}$. We claim that $(g-1)v\in
L(V)\cap V(i+j)+V(i+j+1)$ for all $v\in L^{2k+1}(V)\cap
R^{2(s_i-k)+1}(V)\cap V(i)+V(i+1)$, and that $(g-1)v+V(i+j+1)$
depends only on $v+V(i+j)$, that is, $g_i$ is a well-defined
function $S_{2k+1}^i\to S_1^{i+j}$.

Since $N_j\cap M_{2k-1}$ and $N_{j+1}$ are normal subgroups of $G$
we may write $g=g_1g_2$, where $g_1\in N_j\cap M_{2k-1}$ and
$g_2\in N_{j+1}$. Then $$
g-1=g_1g_2-1=g_1(g_2-1+1)-1=g_1(g_2-1)+(g_1-1).
$$
Suppose first $v\in V(i+1)$. Then from $g\in N_j$ it follows
$(g-1)v\in V(i+j+1)$. Suppose next $v\in L^{2k+1}(V)\cap
R^{2(s_i-k)+1}(V)\cap V(i)$. Then $g_2\in N_{j+1}$ implies
$(g_2-1)v\in V(i+j+1)$, while Theorem \ref{sent-to-L} and the
definition of $N_j$ give $(g_1-1)v\in L(V)\cap
R^{2(s_i-k)+1}(V)\cap V(i+j)$. But $i\in I(j,k)$, so $s_i-k\geq
s_{i+j}$ and therefore Lemma \ref{tower} gives
$$
L(V)\cap R^{2(s_i-k)+1}(V)\cap V(i+j)+V(i+j+1)=L(V)\cap
V(i+j)+V(i+j+1).
$$
Thus $(g_2-1)v\in L(V)\cap V(i+j)+V(i+j+1)$. Our claim now follows
from the above considerations.

We next claim that $g_i$ depends only on the the class
$[g]=g(N_j\cap M_{2k+1})N_{j+1}$ of $g$. For this purpose let
$h\in (N_j\cap M_{2k+1})N_{j+1}$. We may again write $h=h_1h_2$,
where $h_1\in N_j\cap M_{2k+1}$ and $h_2\in N_{j+1}$.  Let $v\in
L^{2k+1}(V)\cap R^{2(s_i-k)+1}(V)\cap V(i)$. As above $(h_2-1)\in
V(i+j+1)$, while the very definition of $M_{2k+1}$ ensures that
$(h_1-1)v=0$. It follows that $(h-1)v\in V(i+j+1)$, thereby
proving our claim.

Since it is clear that $g_i$ is not just a function $S_{2k+1}^i\to
S_1^{i+j}$ but also a linear map, what we have proven so far is
that $Y_i$ is a well-defined function.

We next claim that $Y_i$ is a group homomorphism. Indeed, let
$g,h\in (N_j\cap M_{2k-1})N_{j+1}$. Let $v\in L^{2k+1}(V)\cap
R^{2(s_i-k)+1}(V)\cap V(i)$. Then
$$
(gh-1)v+V(i+j+1)=(h-1)v+(g-1)v+(g-1)(h-1)v+V(i+j+1).
$$
Since $(h-1)v\in L(V)\cap V(i+j)+V(i+j+1)\cap V(i+j)$, and $g\in
N_j$ gives $(g-1)V(i+j)\subseteq V(i+j+1)$, we deduce
$(g-1)(h-1)v\in V(i+j+1)$. It follows that
$$
(gh-1)v+V(i+j+1)=(h-1)v+(g-1)v+V(i+j+1),
$$
thereby proving our claim. We only remaining details to check is
that $Y_i$ commutes with the actions of $G$ and $F$, but this is
straightforward and we can safely omit the details.

{\sc Step II:} $Y_i$ restricted to a subgroup $R_i$ is an
isomorphism.

We again let $i\in I(j,k)$. By Theorem \ref{butterfly} we know
that $X(i,j,2k-1)$ is contained in $(N_j\cap M_{2k-1})N_{j+1}$.
Let $R_i$ be the subgroup of $(N_j\cap M_{2k-1})N_{j+1}/N_j\cap
M_{2k+1})N_{j+1}$ generated by the classes of all elements in
$X(i,j,2k-1)$. The construction of $R_i$ along with Definitions
\ref{bert1} and \ref{doce} show that
\begin{equation}
\label{qi} R_i\subseteq N_{i,j,2k-1}/N_{j,2k+1}.
\end{equation}
Moreover, the definition of $Y_i$ along with Definition
\ref{bert1} show that the images under $Y_i$ of the elements
\begin{equation}
\label{HJ} g^{i,i+j,p,q}_{2k+1,1}N_{j,2k+1}\in R_i,\quad 1\leq
p\leq m_i, 1\leq q\leq m_{i+j},
\end{equation}
form an $F$-basis of $\mathrm{Hom}_F (S^i_{2k+1},S_1^{i+j})$. From
Theorem \ref{irrvi} we know that this space is
$m_im_{i+j}$-dimensional. But the $m_im_{i+j}$ elements (\ref{HJ})
generate $R_i$ as a vector space. It follows that the restriction
of $Y_i$ to $R_i$ is an isomorphism and
\begin{equation}
\label{HJ2} \dim R_i=m_im_{i+j}.
\end{equation}

{\sc Step III:} $Y$ is an isomorphism.

Let
$$Y:N_{j,2k-1}/N_{j,2k+1}\to\underset{i\in
I(j,k)}\oplus\mathrm{Hom}_F (S^i_{2k+1},S_1^{i+j})$$ be the
homomorphism of $FG$-modules induced by the $Y_i$, $i\in I(j,k)$.
The very definitions of $Y_i$ and $N_{i,j,2k-1}$ show that
\begin{equation}
\label{ulti} N_{i,j,2k-1}/N_{j,2k+1}\subseteq\ker Y_{i'},\quad
i\neq i'\in I(j,k).
\end{equation}
We deduce from (\ref{qi}) that
$$
R_i \subseteq\ker Y_{i'},\quad i\neq i'\in I(j,k).
$$
It follows that the image under $Y$ of the product of the
subgroups $R_i$ of $N_{j,2k-1}/N_{j,2k+1}$, as $i$ ranges through
$I(j,k)$, is equal to $\underset{i\in I(j,k)}\oplus\mathrm{Hom}_F
(S^i_{2k+1},S_1^{i+j})$. Since each summand in this space has
dimension $m_im_{i+j}$, the entire space has dimension
$d_{j,2k-1}$. But from Theorem \ref{butterfly} we see that
$$
\dim N_{j,2k-1}/N_{j,2k+1}\leq d_{j,2k-1}.
$$
Since by above $Y$ is an epimorphism, we deduce that $Y$ is an
isomorphism and
\begin{equation}
\label{duy} \dim N_{j,2k-1}/N_{j,2k+1}=d_{j,2k-1}.
\end{equation}

{\sc Step IV:} $N_{j,2k-1}/N_{j,2k+1}$ is internal direct product
of the $R_i$.

By Theorem \ref{butterfly} we know that $N_{j,2k-1}/N_{j,2k+1}$ is
generated as an $F$-vector space by its subspaces $R_i$, $i\in
I(j,k)$. We infer from (\ref{HJ2}) and (\ref{duy}) that as a
vector space $N_{j,2k-1}/N_{j,2k+1}$ is the direct sum of the
$R_i$, therefore as groups we have the following internal direct
product decomposition
$$
N_{j,2k-1}/N_{j,2k+1}=\underset{i\in I(j,k)}\Pi R_i.
$$

{\sc Step V:} $N_{i,j,2k-1}/N_{j,2k+1}=R_i$ for all $i\in I(j,k)$.

Let $i\in I(j,k)$. In view of (\ref{qi}) and (\ref{HJ2}) it
suffices to prove that $\dim N_{i,j,2k-1}/N_{j,2k+1}\leq
m_im_{i+j}$. For this purpose let $P_i$ denote the product of all
$N_{i',j,2k-1}/N_{j,2k+1}$, $i\neq i'\in I(j,k)$. From
(\ref{ulti}) we see that $P_i$ is contained in the kernel of
$Y_i$. This fact and a new application of (\ref{ulti}) yield that
$P_i\cap (N_{i,j,2k-1}/N_{j,2k+1})$ is contained in the kernel of
$Y$. But $Y$ is an isomorphism, so
\begin{equation}
\label{dp} P_i\cap (N_{i,j,2k-1}/N_{j,2k+1})=1,\quad i\in I(j,k).
\end{equation}
But $P_i$ contains all classes of elements in $X(i',j,2k-1)$,
$i\neq i'\in I(j,k)$, so by Theorem \ref{butterfly} the dimension
of the quotient space of $N_{j,2k-1}/N_{j,2k+1}$ by $P_i$ has
dimension at most $m_im_{i+j}$. This and (\ref{dp}) imply
$$
\dim N_{i,j,2k-1}/N_{j,2k+1}\leq m_im_{i+j},
$$
as required. This completes the proof of the theorem.

\begin{thm}
\label{dim-U} $\dim U=\underset{1\leq i<j\leq t}\sum (s_i-s_j+1)
m_im_j.$
\end{thm}

\noindent{\it Proof.} By Theorem \ref{G-action} we have
$$
\mathrm{dim}\,U=\underset{1\leq j}\sum\;\underset{0\leq
k}\sum\;\underset{i\in I(j,k)}\sum m_im_{i+j}=\underset{1\leq
i<j\leq t}\sum (s_i-s_j+1)m_im_j.
$$

\begin{defn} For each $j\geq 1$ let $k(j)$ be the largest integer $k$ such
that $I(j,k)$ is non-empty.
\end{defn}

\begin{thm}
\label{main} There is a canonical $G$-invariant descending central
series for $N/G[V_\infty]$ each of whose factors is naturally an
irreducible $FG$-module, which can be obtained as follows.

We start with $G$-invariant decreasing central series
$$(N_{j,2k-1}/G[V_\infty])_{1\leq j,0\leq k}$$ of $N/G[V_\infty]$
defined in (\ref{nbv}) and refine it by means of the decomposition
(\ref{nice}). The only non-trivial factors thus arising are
\begin{equation}
\label{TY} N_{i,j,2k-1}/N_{j,2k-1}\cong_{FG} \mathrm{Hom}_F
(S^i_{2k+1},S_1^{i+j}),
\end{equation}
where $1\leq j<t-1$, $0\leq k\leq k(j)$, $i\in I(j,k)$, the
dimension of $\mathrm{Hom}_F (S^i_{2k+1},S_1^{i+j})$ is
$m_im_{i+j}$, and the $S^i_{2k+1}$, $S_1^{i+j}$ are amongst the
irreducible constituents of the $FG$-module $V_\infty$ determined
in Theorem \ref{irrvi}.

Each of the factors (\ref{TY}) is an irreducible $FG$-module whose
isomorphism type depends only on $i$ and $j$, and whose
multiplicity in the series is exactly $s_i-s_j+1$. Moreover,
$G[N_{i,j,2k-1}/N_{j,2k-1}]$ contains $N$ and all $E_l$, where
$1\leq l\leq t$, $l\neq i,i+j$, and as a module over
$$G/G[N_{i,j,2k-1}/N_{j,2k-1}]\cong E_{i+j}\times E_i\cong \GL_{m_{i+j}}(F)\times\GL_{m_i}(F)$$
we have
$$
N_{i,j,2k-1}/N_{j,2k-1}\cong M_{m_{i+j}m_i}(F),
$$
where the action is given by
$$
(X,Y)\cdot A=XAY^{-1},\quad X\in GL_{m_{i+j}}(F), A\in
M_{m_{i+j}m_i}(F), Y\in \GL_{m_i}(F).
$$
\end{thm}

\noindent{\it Proof.} By Theorem \ref{G-action} the only
non-trivial factors of the series $(N_{j,2k-1}/G[V_\infty])_{1\leq
j,0\leq k}$ are of the form $\mathrm{Hom}_F
(S^i_{2k+1},S_1^{i+j})$, where $1\leq j<t-1$, $0\leq k\leq k(j)$,
$i\in I(j,k)$. From Theorem \ref{irrvi} we know that $N$ acts
trivially on all this factors, that $E_l$ also acts trivially on
them if $l\neq i,i+j$, and that $E_{i+j}\times E_i$ acts
irreducibly as indicated. Since such factor appears as many times
as $k$ is between $0$ and $s_i-s_{i+j}$, the result follows.

\section{The split extension $G[V_\infty]$ of $G[V_\infty]\cap G[{}^\infty
V/V_\infty]$}

\begin{thm}
\label{bigsplit} The canonical restriction map $G[V_\infty]\to
G({}^\infty V/V_\infty)$ is a split group epimorphism with kernel
$G[V_\infty]\cap G[{}^\infty V/V_\infty]$ and complement
$G(V_{\even})\times G(V_{\ndeg})$. Moreover, $U$ normalizes $E$,
so that
\begin{equation}
\label{26-1} G\{V_{\odd}^\dag\}\cap G(V_{\odd})=U\rtimes E,
\end{equation}
\begin{equation}
\label{26-2} G=\GI\rtimes (U\rtimes E),
\end{equation}
and $$ G=(G[V_\infty]\cap G[{}^\infty V/V_\infty])\rtimes
\left(G(V_{\even})\times G(V_{\ndeg})\times (U\rtimes E)\right),
$$
where the structure of $U$ as a group under the action of $E$ is
described in section \ref{nginfty}.
\end{thm}

\noindent{\it Proof.} We know from \cite{DS1} that $G({}^\infty
V/V_\infty)$ preserves the even and non-degenerate parts of
${}^\infty V/V_\infty$. It follows that the restriction map
$$
G(V_{\even})\times G(V_{\ndeg})\hookrightarrow G[V_\infty] \to
G({}^\infty V/V_\infty)
$$
is a group isomorphism. Since
$$
G(V_{\even})\times G(V_{\ndeg})\cap (G[V_\infty]\cap G[{}^\infty
V/V_\infty])=1,
$$
we infer
$$
G[V_\infty]=(G[V_\infty]\cap G[{}^\infty V/V_\infty])\rtimes
(G(V_{\even})\times G(V_{\ndeg})).
$$
The very characterizations of $E$ and $U$ given in Theorems
\ref{E} and \ref{U} show that $E$ normalizes $U$, and both groups
are contained in $G\{V_{\odd}^\dag\}\cap G(V_{\odd})$. Moreover,
it is obvious that $G(V_{\even})\times G(V_{\ndeg})$ commutes
elementwise with $G(V_{\odd})$. Furthermore, from Lemma \ref{incl}
we deduce
$$
G[V_\infty]\cap G\{V_{\odd}^\dag\}\cap G(V_{\odd})=1.
$$
The desired conclusion now follows from Theorems \ref{E} and
\ref{U}.

\section{A criterion applicable to bilinear spaces of types E and I}
\label{go}

We next derive a criterion that yields the structure of the
isometry group of a bilinear space of type  E, i.e. the space is
equal to its even part, or type I, in C. Riehm's notation.

\begin{notn} If $Y$ is an $F$-vector space and
$u\in\End(Y)$ then $C_{\GL(Y)}(u)$ denotes the centralizer of $u$
in $\GL(Y)$.
\end{notn}

\begin{notn} If $Y$ and $Z$ are $F$-vector spaces then $\Bil(Y,Z)$
denotes the $F$-vector space of all bilinear forms $Y\times Z\to
F$. We say that $\phi\in\Bil(Y,Z)$ is non-degenerate if its left
and right radicals are equal to (0).
\end{notn}

\begin{thm}
\label{YZ} Let $(W,\phi)$ be a bilinear space. Suppose there
exists $G(W,\phi)$-invariant totally isotropic subspaces $Y$ and
$Z$ of $W$ such that $W=Y\oplus Z$ and $\phi|_{Z\times Y}$ is
non-degenerate. Then

(1) There exists a unique $u\in \End_F(Y)$ such that
\begin{equation}
\label{uyz} \phi(y,z)=\phi(z,uy),\quad y\in Y,z\in Z.
\end{equation}

(2) If $g\in G(W,\phi)$ then $g|_Y\in C_{\GL(Y)}(u)$.

(3) The canonical restriction map $\rho:G(W,\phi)\to
C_{\GL(Y)}(u)$, given by $g\mapsto g|_Y$, is a group isomorphism.
\end{thm}

\noindent{\it Proof.} Consider the linear map $A:\End_F(Y)\to
\Bil(Y,Z)$, given by $u\mapsto \phi_u$, where
$$
\phi_u(y,z)=\phi(z,uy), \quad y\in Y,z\in Z.
$$
Since the right radical of $\phi|_{Z\times Y}$ is $(0)$, it
follows that $A$ is a monomorphism. As the left radical of
$\phi|_{Z\times Y}$ is also $(0)$, we infer $\dim Y=\dim Z$,
whence $\dim \End_F(Y)=\dim \Bil(Y,Z)$, so $A$ is an isomorphism.
In particular, there exists a unique $u\in \End_F(Y)$ such that
$A(u)=\phi|_{Y\times Z}$.

Let $g\in G(W,\phi)$. For $y\in Y$ and $z\in Z$, since both $Y$
and $Z$ are $G(W,\phi)$-invariant, (\ref{uyz}) gives
$$
\phi(z,guy)=\phi(g^{-1}z,uy)=\phi(y,g^{-1}z)=\phi(gy,z)=\phi(z,ugy).
$$
As the right radical of $\phi|_{Z\times Y}$ is $(0)$, we deduce
$g|_Y\in C_{\GL(Y)}(u)$.

Let $g\in\ker \rho$. For $y\in Y$ and $z\in Z$ we have
$$
\phi(gz,y)=\phi(gz,gy)=\phi(z,y).
$$
As the left radical of $\phi|_{Z\times Y}$ is $(0)$, we obtain
$g|_Z=1_Z$. But $W=Y\oplus Z$, so $g=1$. This proves that $\rho$
is injective.

Let $b\in C_{\GL(Y)}(u)$. Consider the linear map $\End_F(Z)\to
\Bil(Z,Y)$, given by $c\mapsto \phi^c$, where
$$
\phi^c(z,y)=\phi(cz,y), \quad y\in Y,z\in Z.
$$
As above, this is an isomorphism. In particular, there exists a
unique $c\in \End_F(Z)$ such that
$$
\phi(cz,y)=\phi(z,b^{-1}y), \quad y\in Y,z\in Z.
$$
As $b\in \GL(V)$ and left radical of $\phi|_{Z\times Y}$ is $(0)$,
we infer that $c\in \GL(Z)$. We may re-write the above equation in
the form
\begin{equation}
\label{rew} \phi(cz,by)=\phi(z,y), \quad y\in Y,z\in Z.
\end{equation}
Let $g=b\oplus c\in \GL(W)$. Let $y_1,y_2\in Y$ and $z_1,z_2\in
Z$. Since $Y$ and $Z$ are totally isotropic, (\ref{uyz}) and
(\ref{rew}) along with $b\in C_{\GL(Y)}(u)$ give
$$
\begin{aligned}
\phi(g(y_1+z_1),g(y_2+z_2)) &=\phi(by_1+cz_1,by_2+cz_2)=
\phi(by_1,cz_2)+ \phi(cz_1,by_2)\\
&=\phi(cz_2,uby_1)+ \phi(z_1,y_2)=\phi(cz_2,buy_1)+
\phi(z_1,y_2)\\
&=\phi(z_2,uy_1)+ \phi(z_1,y_2)=\phi(y_1,z_2)+ \phi(z_1,y_2)\\
&=\phi(y_1+z_1,z_2)+ \phi(y_1+z_1,y_2)=\phi(y_1+z_1,y_2+z_2).
\end{aligned}
$$
Therefore $g\in G(W,\phi)$. By construction $\rho(g)=b$, so $\rho$
is an epimorphism, thus completing the proof.

\begin{note}
\label{dual}
 Observe that $Z^*$, the dual of $Z$, is an $FG(W,\phi)$-module
in the usual way. Moreover, $Y$ and $Z^*$ are isomorphic, as
$FG(W,\phi)$-modules, via the map $y\mapsto \phi(-,y)$.
\end{note}

\subsection{Structure of $G(V_\even)$}

\begin{thm}
\label{ee} Suppose that $V=V_\even$. Then

(1) For uniquely determined positive integers $n_i$ and $r_i$ we
have an equivalence of bilinear spaces
$$
V\cong \underset{1\leq i\leq d}\perp n_i N_{2r_i}.
$$

(2) There exists a unique $u\in \End_F(L_\infty(V))$ such that
$$
\vf(l,r)=\vf(r,ul),\quad l\in L_\infty(V),r\in R_\infty(V).
$$
The endomorphism $u$ is nilpotent, with elementary divisors
$t^{r_1},...,t^{r_d}$ and multiplicities $n_1,...,n_d$.

(3) The canonical restriction map $\rho:G(V)\to
C_{\GL(L_\infty(V))}(u)$, given by $g\mapsto g|_{L_\infty(V)}$, is
a group isomorphism.
\end{thm}

\noindent{\it Proof.} The first assertion follows from Theorem
\ref{gabriel}. By means of Lemmas \ref{block2} and \ref{LandR} we
deduce that the subspaces $L_\infty(V)$ and $R_\infty(V)$ of $V$
satisfy the hypotheses of Theorem \ref{YZ}. This theorem yields
all remaining assertions, except for the similarity type of $u$.
By hypothesis there is a basis of $V$ relative to which the matrix
of $\vf$ is equal to $\underset{1\leq i\leq d}\oplus n_i
J_{2r_i}(0)$. A suitable rearrangement of this basis which puts
first all basis vectors of $R_\infty(V)$ and second all basis
vectors of $L_\infty(V)$ yields a new basis relative to which the
matrix of $\vf$ is equal to
$$\left(\begin{array}{cc}0&1\\J&0\end{array}\right),
$$
where $J=\underset{1\leq i\leq d}\oplus n_i J_{r_i}(0)$. Since $J$
is the matrix of $u$ in the above basis of $L_\infty(V)$, the
similarity type of $u$ is as given.

\subsection{Irreducible constituents of the $FG(V_\even)$-module
$V_\even$}

We know from Note \ref{dual} that $R_\infty(V)^*\cong
L_\infty(V)$, as $FG$-modules. In view this, the decomposition
$V=L_\infty(V)\oplus R_\infty(V)$ and Theorem \ref{ee}, it
suffices to restrict ourselves to the classical case of finding
the irreducible constituents of $L_\infty(V)$ as an
$FC_{\GL(L_\infty(V))}(u)$-module. This is well known and will be
omitted.

\section{Structure of $G(V_{\ndeg})$}
\label{to}

We assume here that $V=V_\ndeg$. Note that $\Bil(V)$ is a natural
right $\End(V)$-module via
$$
(\phi\cdot u)(x,y)=\phi(x,uy),\quad \phi\in\Bil(V),u\in
\End(V),x,y\in V.
$$
For a fixed $\phi\in\Bil(V)$ the map $\End(V)\to\Bil(V)$ given by
$u\to \phi\cdot u$ is a linear isomorphism if and only if $\phi$
is non-degenerate, in which case $u$ is invertible if and only if
$\phi\cdot u$ is non-degenerate. In this case, given any
$\psi\in\Bil(V)$ we shall write $u_{\phi,\psi}$ for the unique
$u\in \End(V)$ such that $\psi(x,y)=\phi(x,uy)$.

Since $\vf$ is non-degenerate, we may use it to represent any
bilinear form, in particular $\vf'$. We write $\s=u_{\vf,\vf'}$
for the {\em asymmetry} of $\vf$, i.e. the element of $\GL(V)$
satisfying
\[ \vf'(x,y)=\vf(x,\sigma(y)),\quad x,y\in V. \]
This linear operator measures how far is $\vf$ from being
symmetric. We have
$$\vf(x,y)=\vf(y,\sigma(x))=\vf(\sigma(x),\sigma(y))\quad x,y\in
V,
$$
so that $\sigma\in G$. In fact, it is easy to see that $\sigma$
belongs to the center $Z(G)$ of $G$.

Let $F[t]$ denote the polynomial algebra in one variable $t$ over
$F$. We view $V$ as an $F[t]$-module via $\sigma$. For $0\neq q\in
F[t]$, consider the adjoint polynomial $q^*\in F[t]$, defined by
$$q^*(t)=t^{\mathrm{deg}\,q}q(1/t).$$

The minimal polynomial of $\sigma$ will be denoted by $p_\sigma\in
F[t]$. Let $\P$ stand for the set of all monic irreducible
polynomials in $F[t]$ dividing $p_\sigma$. For $p\in \P$ let $V_p$
denote the primary component of $\sigma$ associated to $p$. Since
$\sigma\in Z(G)$, each primary component is $G$-invariant. We
consider the subsets of $\P$:
$$
\P_1=\{p\in \P\,|\, p^*\neq \pm p\}\text{ and }\P_2=\{p\in
\P\,|\,p^*=\pm p\}.$$ We construct a subset $\P_1'$ of $\P_1$ by
selecting one element out of each set $\P_1\cap \{\pm p,\pm
p^*\}$, as $p$ ranges through $\P_1$. It follows at once from
\cite{CR} that
\begin{equation}
\label{dw} G(V)\cong \left(\underset{p\in \P_1'}\Pi G(V_p\oplus
V_{p^*})\right)\Pi\left(\underset{p\in \P_2}\Pi G(V_p)\right)
\end{equation}
Thus the study of $G$ reduces to two cases:

{\sc Case I:} $V=V_p\oplus V_{p^*}$, $p^*\neq \pm p$.

{\sc Case II:} $V=V_p$, $p^*=\pm p$.

We break II up into two cases:

{\sc Case II}a: $\mathrm{deg}\,p>1$ or $\mathrm{char}\,F\neq 2$.

{\sc Case II}b: $\mathrm{deg}\, p=1$ and $\mathrm{char}\,F=2$.

\subsection{Case I}

We assume here that $p$ is a monic irreducible polynomial in
$F[t]$ dividing $p_\sigma$ such that $p^*\neq \pm p$ and
$V=V_p\oplus V_{p^*}$. In particular, $(p,p^*)=1$. As shown in
\cite{CR} the $G$-invariant $F[t]$-submodules $V_p$ and $V_{p^*}$
of $V$ are totally isotropic. In view of Theorem \ref{YZ}, we have
the following result.

\begin{thm} \label{Ha} The restriction map
$\rho: G\to C_{\GL(V_p)}(\sigma|_{V_p})$ is an isomorphism.
\end{thm}

Note that when $F$ is algebraically closed $p=t-\lambda$ for some
$\lambda\in F$ different from $1$ and $-1$. In this case then $G$
becomes isomorphic to the centralizer of a nilpotent element in
the general linear group (as adding a scalar operator does not
change the centralizer).

\subsection{Case IIa}

We assume here that $p$ is a monic irreducible polynomial in
$F[t]$ dividing $p_\sigma$ such that $p^*=\pm p$ and $V=V_p$. We
further assume that $\mathrm{deg}\,p>1$ or $\mathrm{char}\,F\neq
2$. The {\em symmetric} and {\em alternating} parts of $\vf$ are
defined by $\vf^\pm=\vf\pm \vf'$. Clearly $G(\vf)\subseteq
G(\vf^+)\cap G(\vf^-)$, with equality if $\mathrm{char}\,F\neq 2$.

\begin{lem}
\label{nsym} $\vf^\pm$ is non-degenerate if and only if $p_\s(\mp
1)\neq 0$.
\end{lem}

\noindent{\it Proof.} This follows from the identity
$$
\vf^\pm(x,y)=\vf(x,(1\pm \sigma)y),\quad x,y\in V.
$$

If $\vf^+$ is non-degenerate, we write $\s^{+-}=u_{\vf^+,\vf^-}$
and $\s^+=u_{\vf^+,\vf}$; moreover, we denote the isometry group
of $\vf^+$ by $\O(\vf^+)$ and the associated Lie algebra by
$\o(\vf^+)$. If $\vf^-$ is non-degenerate, we write
$\s^{-+}=u_{\vf^-,\vf^+}$ and $\s^-=u_{\vf^-,\vf}$; moreover, we
denote the isometry group of $\vf^-$ by $\Sp(\vf^-)$ and the
associated Lie algebra by $\sp(\vf^+)$.

\begin{lem}
 \label{lie0}
 If $\vf^+$ is non-degenerate then $\s^{+-}\in
\o(\vf^+)$. If  $\vf^-$ is non-degenerate then $\s^{-+}\in
\sp(\vf^-)$.
\end{lem}

\noindent{\it Proof.} We have
$$
\vf^+(\s^{+-}x,y)+\vf^+(x,\s^{+-}y)=\vf^+(y,\s^{+-}x)+\vf^-(x,y)=\vf^-(y,x)+\vf^-(x,y)=0,\quad
x,y\in V,
$$
thereby proving the first assertion. The second is proved
similarly.

\begin{prop}
 \label{lie1}
Suppose $\vf^+$ is non-degenerate. Then
$$
G(\vf)=C_{\O(\vf^+)}(\s^+),
$$
and if $\mathrm{char}\,F\neq 2$ then
$$
G(\vf)=C_{\O(\vf^+)}(\s^{+-}).
$$
\end{prop}

\noindent{\it Proof.} Let $a\in \GL(V)$. We have
$$a\in G(\vf)$$
if and only if
$$
\vf(ax,ay)=\vf(x,y),\quad x,y\in V
$$
if and only if
$$
\vf^+(ax,\s^+ ay)=\vf^+(x, \s^+ y),\quad x,y\in V
$$
if and only if
$$
\vf^+(x,a^{-1}\s^+ ay)=\vf^+(x, \s^+ y),\; x,y\in V\text{ and
}a\in \O(\vf^+)
$$
if and only if
$$
a\in C_{\O(\vf^+)}(\s^+).
$$
This proves the first assertion. As for the second, if
$\mathrm{char}\,F\neq 2$ then
$$a\in G(\vf)$$
if and only if
$$a\in \O(\vf^+)\text{ and }a\in G(\vf^-)$$
if and only if
$$a\in \O(\vf^+)\text{ and }\vf^-(ax,ay)=\vf^-(x,y),\quad x,y\in V$$
if and only if
$$a\in \O(\vf^+)\text{ and }\vf^+(ax,\sigma^{+-}ay)=\vf^+(x,\sigma^{+-}y),\quad x,y\in
V$$ \ if and only if
$$a\in \O(\vf^+)\text{ and }\vf^+(x,a^{-1}\sigma^{+-}ay)=\vf^+(x,\sigma^{+-}y),\quad x,y\in
V$$ if and only if
$$a\in C_{\O(\vf^+)}(\s^{+-}).$$

\begin{prop}
 \label{lie2}
Suppose $\vf^-$ is non-degenerate. Then
$$
G(\vf)=C_{\Sp(\vf^-)}(\s^-),
$$
and if $\mathrm{char}\,F\neq 2$ then
$$
G(\vf)=C_{\Sp(\vf^-)}(\s^{-+}).
$$
\end{prop}

\noindent{\it Proof.} This is similar to the above proof, mutatis
mutandi.

We know from \cite{CR} that if $\mathrm{deg}\,p>1$ then
$\mathrm{deg}\,p$ is even and $p=p^*$, while it is obvious that if
$\mathrm{deg}\,p=1$ then $p=t\pm 1$.

\begin{thm}
\label{gral}
 (i) If $\mathrm{deg}\,p>1$ then $\vf^\pm$ is non-degenerate, $\s^{+-}\in
\o(\vf^+)$, $\s^{-+}\in \sp(\vf^-)$ and
$$
G=C_{\O(\vf^+)}(\s^+)=C_{\Sp(\vf^-)}(\s^-).
$$
Moreover, if $\mathrm{char}\,F\neq 2$ then
$$
G=C_{\O(\vf^+)}(\s^{+-})=C_{\Sp(\vf^-)}(\s^{-+}).
$$

(ii) If $p=t-1$ and $\mathrm{char}\,F\neq 2$ then $\vf^+$ is
non-degenerate, $\s^{+-}\in \o(\vf^+)$ and
$$
G=C_{\O(\vf^+)}(\s^+)=C_{\O(\vf^+)}(\s^{+-}).
$$

(iii) If $p=t+1$ and $\mathrm{char}\,F\neq 2$ then $\vf^-$ is
non-degenerate, $\s^{-+}\in \sp(\vf^-)$ and
$$
G=C_{\Sp(\vf^-)}(\s^-)=C_{\Sp(\vf^-)}(\s^{-+}).
$$
\end{thm}

\noindent{\it Proof.} This follows from Lemmas \ref{nsym} and
\ref{lie0}, and Propositions \ref{lie1} and \ref{lie2}.

For the remainder of this subsection we suppose that $F$ is
algebraically closed of characteristic not 2. Then $p=t\pm 1$.

For convenience, we define $n$-by-$n$ matrices $H_n(\lambda)$ and
$\Gamma_n$ by
\[ H_n(\lambda)=\left(\begin{array}{cc}0&I_m\\ J_m(\lambda)
&0\end{array}\right), \quad n=2m, \lambda\in F, \] and
\[ \Gamma_n=\left(\begin{array}{crrccccccc}
0&0&0&0&\cdots &0&0&(-1)^{n-1}\\
0&0&0&0&\cdots &0&(-1)^{n-2}&(-1)^{n-2}\\
\vdots &&&&&&\\
0&-1&-1&0&\cdots &0&0&0 \\
1&1&0&0&\cdots &0&0&0\end{array}\right). \]

We refer the reader to \cite{HS}, and for an older version,
\cite{CW}, for a proof of the Canonical Form Theorem for bilinear
forms.

\begin{thm} \label{class}
(a) Any $\phi\in\Bil(V)$ admits an orthogonal direct decomposition
\[ \phi=\phi_1\perp\phi_2\perp\cdots\perp\phi_k, \]
where the $\phi_i$'s are indecomposable bilinear forms which are
unique up to equivalence and permutation.

(b) If $\phi\in\Bil(V)$ is indecomposable then, with respect to a
suitable basis of $V$, the matrix of $\phi$ is one of the
following:
\begin{itemize}
\item[(i)] $H_n(\lambda)$, $n=2m$, $\lambda\ne(-1)^{m+1}$;
\item[(ii)] $\Gamma_n$, $n\ge1$; \item[(iii)] $J_n(0)$, $n=2m+1$.
\end{itemize}

(c) The matrices listed in part (b) are pairwise non-congruent
except for the fact that $H_n(\lambda)$ and $H_n(\lambda^{-1})$
are congruent when $\lambda\ne0,\pm1$.
\end{thm}

We mention that, when $n=2m$ is even, $H_n(0)$ is congruent to
$J_n(0)$.

\begin{thm} \label{sv+1}
If $p=t-1$ then $\vf^+$ is non-degenerate, $\s^{+-}$ belongs to
$\so(\vf^+)$, $G=C_{\O(\vf^+)}(\s^{+-})$, and the linear operators
$\sigma-1$ and $\s^{+-}$ are nilpotent and similar to each other
(i.e., they have the same elementary divisors).
\end{thm}

\noindent{\it Proof.} In view of Theorem \ref{gral} we are reduced
to show the last assertion. It suffices to verify this assertion
for indecomposable $\vf$. There are two cases to consider. The
matrix of $\vf$ will be denoted by $A_\vf$.

First, the matrix of $\vf$ is $H_n(1)$ where $n=2m$ and $m$ is
even. Then the matrix of $u$ is $-A_{\vf^+}^{-1} A_{\vf^-}$. An
easy computation shows that both $\sigma-1$ and $\s^{+-}$ have
elementary divisors $t^m$ and $t^m$.

Second, the matrix of $\vf$ is $\Gamma_n$ and $n$ is odd. In that
case the matrix $A_{\vf^+}$ is involutory and a simple computation
shows that the matrix of $\s^{+-}$ is equal to $-J_n(0)'$. Hence
both $\sigma-1$ and $\s^{+-}$ have only one elementary divisor,
namely $t^n$.

\begin{thm} \label{sv-1}
If $p=t+1$ then $\vf^-$ is non-degenerate, $\s^{-+}$ belongs to
$\sp(\vf^-)$, $G=C_{\Sp(\vf^-)}(\s^{-+})$, and the linear
operators $\sigma+1$ and $\s^{-+}$ are nilpotent and similar to
each other (i.e., they have the same elementary divisors).
\end{thm}

\noindent{\it Proof.} This proof is similar to the one above.

\subsection{Case IIb}

We have not been able to make progress on this case.

\section{The 2-step nilpotent group $G[V_\infty]\cap G[{}^\infty
V/V_\infty]$}

The divide the study of $\GIK$ into two cases, namely that of
$G[\VI]$ and $\GIK/G[\VI]$.

\subsection{Basic facts about $\GIK$}


\begin{lem}
\label{step2} $G[\VI]$ is contained in the center of $\GIK$.
\end{lem}

\noindent{\it Proof.} From Lemma \ref{step1} we know that $G[\VI]$
is contained in $G[V_\infty]\cap G[V/V_\infty]$, which is in turn
contained in $\GIK$. Let $g\in G[\VI]$ and $h\in \GIK$. By Lemma
\ref{step1} we have $(g-1)V\subseteq V_\infty$. Since $\GIK$ is
the identity on $V_\infty$, we infer $(h-1)(g-1)V=(0)$. Moreover,
$(h-1)V\subseteq {}^\infty V$ by Lemma \ref{incl}, whence
$(g-1)(h-1)V=0$ by the very definition of $G[\VI]$. It follows
that
\begin{equation}
\label{gh}
 hg=h+g-1=gh,
\end{equation}
as required.

\begin{lem}
\label{step3} $\GIK$ is nilpotent of class $\leq 2$.
\end{lem}

\noindent{\it Proof.} By Lemma \ref{step2} it suffices to show
that $\GIK/G[\VI]$ is abelian. Let $g,h\in\GIK$ and let $v\in
{}^\infty V$. Since $gv-v\in V_\infty$, by the definition of
$\GIK$, and $h$ is the identity on $V_\infty$, we have
$$h(g(v))=h(g(v)-v+v)=gv-v+hv.
$$
For the same reasons as above
$$
(h^{-1}g^{-1}hg)(v)=h^{-1}g^{-1}(gv+hv-v)=h^{-1}(v+hv-v)=h^{-1}hv=v.
$$
Thus $[h,g]\in G[\VI]$, as required.

\begin{lem}
\label{newstep} $G[V_\infty]\cap G[V/V_\infty]$ is an abelian
unipotent normal subgroup of $G$. In fact, if $g,h\in
G[V_\infty]\cap G[V/V_\infty]$ then
\begin{equation}
\label{dos} (h-1)(g-1)=0.
\end{equation}
\end{lem}

\noindent{\it Proof.} Since $V_\infty$ is $G$-invariant, it
follows that $G[V_\infty]\cap G[V/V_\infty]$ is a normal subgroup
of $G$. If if $g,h\in G[V_\infty]\cap G[V/V_\infty]$ then
$(g-1)V\subseteq V_\infty$, so $(h-1)(g-1)V=0$. This completes the
proof.

\begin{lem}
\label{step4} $\GIK$ is a unipotent normal subgroup of $G$. In
fact, if $g,h,k\in \GIK$ then $(k-1)(h-1)(g-1)=(0)$.
\end{lem}

\noindent{\it Proof.} Since ${}^\infty V$ and $V_\infty$ are
$G$-invariant, it follows that $\GIK$ is a normal subgroup of $G$.
Let $g,h,k\in \GIK$. By Lemma \ref{incl} we have $(g-1)V\subseteq
{}^\infty V$. From the very definition of $\GIK$ we obtain
$(h-1)(g-1)V\subseteq V_\infty$ and a fortiori
$(k-1)(h-1)(g-1)V=(0)$. This completes the proof.

\subsection{Structure of the $FG$-module $G[\VI]\cap G[V/\R(V)]$}

The divide the study of $G[\VI]$ into two cases, namely that of
$G[\VI]\cap G[V/\R(V)]$ and $G[\VI]/G[\VI]\cap G[V/\R(V)]$.

\begin{thm}
\label{h1} The map $G[{}^\infty V]\to \End_F(V)$ given by
\begin{equation}
\label{h} g\mapsto g-1
\end{equation}
is a group monomorphism whose image is an $F$-vector subspace of
$\End_F(V)$. By transferring this $F$-vector space structure to
$G[V_\infty]\cap G[V/V_\infty]$ the map (\ref{h}) becomes an
$FG$-module monomorphism. The map (\ref{h}) induces a monomorphism
of $FG$-modules
$$
G[\VI]\to \mathrm{Hom}_F(V/{}^\infty V,V_\infty),
$$
and hence a monomorphism of $F$-vector spaces
\begin{equation}
\label{h6} G[\VI]\to \mathrm{Hom}_F(V_\odd^\dag,V_\infty),
\end{equation}
namely by means of $g\mapsto (g-1)\vert_{V_\odd^\dag}$.
\end{thm}

\noindent{\it Proof.} The identity of (\ref{dos}) shows that
(\ref{h}) is a group homomorphism, which is clearly injective and
preserves the action of $G$. Suppose $g\in G[{}^\infty V]$ and
$k\in F$. Then $k(g-1)+1$ is a linear automorphism of $V$ which
fixes ${}^\infty V$ pointwise, acts trivially on $V/V_\infty$ and
preserves the orthogonality of the generators $e_{2k}^{i,p}$ of
$V_\odd^\dag$. It follows that $k(g-1)+1\in G[{}^\infty V]$, so
the image of (\ref{h}) is a subspace of $\End_F(V)$. By Lemma
\ref{step1} we know that (\ref{h}) maps $G[\VI]$ into
$\mathrm{Hom}_F(V,V_\infty)$, and the very definition of $G[\VI]$
yields an induced $FG$-monomorphism $G[\VI]\to
\mathrm{Hom}_F(V/{}^\infty V,V_\infty)$. Since $V_\odd^\dag$
complements ${}^\infty V$ in $V$, the last assertion follows.

\begin{lem}
\label{equal} $G[\VI]\cap G[V/L(V)]=G[\VI]\cap G[V/\R(V)]$.
\end{lem}

\noindent{\it Proof.} By definition the right hand side is
contained in the left hand side. Let $g\in G[\VI]\cap G[V/L(V)]$.
We wish to show that $(g-1)V\subseteq \R(V)$. Since $g$ is the
identity on ${}^\infty V$, it suffices to prove
$(g-1)V_\odd^\dag\subseteq \R(V)$. By assumption
$(g-1)V_\odd^\dag\subseteq L(V)$, so are reduced to demonstrate
$(g-1)V_\odd^\dag\subseteq R(V)$. By Lemma \ref{step1} we have
$(g-1)V_\odd^\dag\subseteq V_\infty$, which leave only the
identity $\la V_\odd^\dag,(g-1)V_\odd^\dag\ra=0$ to be shown.
Well, if $v,w\in V_\odd^\dag$ then $gw-w\in L(V)$, so
$$
0=\la w,v\ra=\la gw,gv\ra=\la gw,gv\ra=\la (gw-w)+w,gv\ra=\la
w,gv\ra=\la w,gv-v\ra,
$$
as required.

\begin{defn} Let $\Bil(V,\VI)$ be the $FG$-submodule of $\Bil(V)$
consisting of all bilinear forms whose radical contains $\VI$.
Thus $\Bil(V,\VI)$ and $\Bil(V/\VI)$ are isomorphic as
$FG$-modules.
\end{defn}

\begin{thm}
\label{bil} The map $G[\VI]\to\Bil(V,\VI)$ given by $g\mapsto
\vf_g$, where
$$
\vf_g(v,w)=\vf((g-1)v,w)=\langle(g-1)v,w\rangle,\quad v,w\in V,
$$
is an $FG$-module homomorphism, inducing an $FG$-module
homomorphism $g\mapsto \widehat{\vf_g}$ from $G[\VI]$ to
$\Bil(V/\VI)$. The kernel of both maps is equal to $G[\VI]\cap
G[V/\R(V)]$.
\end{thm}

\noindent{\it Proof.} The fact that $g\mapsto \vf_g$ is a
homomorphism of $FG$-modules is easily verified. By the very
definition of this map its kernel is equal to $G[\VI]\cap
G[V/L(V)]$, which equals $G[\VI]\cap G[V/\R(V)]$ by Lemma
\ref{equal}.

\begin{notn} The image of $G[\VI]$ under the above $FG$-homomorphism
$G[\VI]\to\Bil(V/\VI)$ will be denoted by $S$.
\end{notn}

\begin{lem}
\label{isomo} The restriction of (\ref{h}) to $G[\VI]\cap
G[V/\R(V)]$ yields an isomorphism of $FG$-modules
\begin{equation}
\label{noname} G[\VI]\cap G[V/\R(V)]\to \mathrm{Hom}_F(V/{}^\infty
V,\R(V)).
\end{equation}
\end{lem}

\noindent{\it Proof.} All maps of the form $1_{{}^\infty V}\oplus
(1_{V_\odd^\dag}+f)$, where $f\in
\mathrm{Hom}_F(V_\odd^\dag,\R(V))$, belong to $G[\VI]\cap
G[V/\R(V)]$, thereby proving that (\ref{noname}) is an
epimorphism. The rest follows from Theorem \ref{h1}.

\begin{thm}
\label{onecase} $G[\VI]\cap G[V/\R(V)]$ is an $FG$-module of
dimension $(\dim V/\VI)\times (\dim \R(V))$. If $\R(V)\neq 0$ and
$V/\VI\neq 0$ its irreducible constituents are of the form
$\mathrm{Hom}_F(Q_{2k}^i,\R(V))$, where the $Q_{2k}^i$ are the
irreducible constituents of the $FG$-module $V/{}^\infty V$
described in Theorem \ref{irrvI}. Each has dimension $m_im_t$,
$1\leq i<t$, and multiplicity $s_i$ with stabilizer
$S_i=N\rtimes\underset{l\neq i,t}\Pi E_l$. As a module for
$G/S_i\cong E_i\times E_t\cong \GL_{m_i}(F)\times \GL_{m_t}(F)$,
$\mathrm{Hom}_F(Q_{2k}^i,\R(V))$ is isomorphic to $M_{m_tm_i}(F)$,
where $(X,Y)\in \GL_{m_i}(F)\times \GL_{m_t}(F)$ acts on $A\in
M_{m_tm_i}(F)$ by $(X,Y)\cdot A= YAX^{-1}$.
\end{thm}

\noindent{\it Proof.} This follows easily from Lemma \ref{isomo}
and Theorem \ref{irrvI}.

We next wish to determine the structure of the the remaining part
of $G[\VI]$, namely
$$
G[\VI]/G[\VI]\cap G[V/\R(V)]\cong_{FG} S.
$$
We digress to record some basic facts from Linear Algebra which
will be required for a complete understanding of the structure of
$S$.

\subsection{$\GL_m(F)$ acting by congruence on $M_m(F)$}

\begin{defn} Let $m\geq 1$. Denote by $S_m(F)$ and $A_m(F)$ the set
of all $m\times m$ symmetric and alternating matrices over $F$,
respectively.
\end{defn}

\begin{thm}
\label{symme}
 Let $m\geq 1$. The irreducible constituents of $M_m(F)$ as a
module for $\GL_m(F)$ over $F$, acting by congruence are as
follows.

(1) If $\mathrm{char}\,F\neq 2$ then
$$
M_m(F)=S_m(F)\oplus A_m(F),
$$
where both summands are irreducible if $m>1$, while
$M_1(F)=S_1(F)$ is irreducible.

(2) If $\mathrm{char}\,F=2$ each factor of the
$\GL_m(F)$-invariant series
$$
0\subseteq A_m(F)\subseteq S_m(F)\subseteq M_m(F)
$$
is irreducible, except when $m=1$ in which case $M_1(F)=S_1(F)$ is
irreducible.

(3) $M_m(F)/S_m(F)$ is isomorphic to $A_m(F)$.

\end{thm}

\noindent{\it Proof.} Consider the homomorphism of
$\GL_m(F)$-modules $M_m(F)\to A_m(F)$, given by $A\mapsto A-A'$.
Since its kernel is $S_m(F)$, a dimension comparison shows that
its image is $A_m(F)$. Therefore $M_m(F)/S_m(F)\cong A_m(F)$. In
view of this isomorphism we may assume throughout that $m>1$, and
we are reduced to show that $A_m(F)$ and $S_m(F)/A_m(F)$ are
irreducible, in the later case when $\mathrm{char}\,F=2$.

{\sc Step I:} $A_m(F)$ is irreducible. If $m=2,3$ there is a
single non-zero $\GL_m(F)$-orbit in $A_m(F)$, which is then
irreducible. Suppose $m>3$. Let $0\neq M$ be any
$F\GL_m(F)$-submodule of $A_m(F)$. To see that $M=A_m(F)$ it
suffices to show that $M$ contains a matrix of rank 2. It is
well-known that $M$ contains a matrix, say $A$, which is the
direct sum of at least one block of the form $\left(\begin{matrix}
0 & 1\\
-1 & 0\end{matrix}\right)$ plus zero blocks. We may assume that
$A$ has at least two non-zero blocks. Choose $B$ in $A_m(F)$ whose
only nonzero entries are in positions $(2,3)$ and $(3,2)$. Then
$A+B$ has the same rank as $A$ and so $A+B$ is in $M$. Hence $B$
is in $M$, as required.

{\sc Step II:} $S_m(F)/A_m(F)$ is irreducible if
$\mathrm{char}\,F=2$. Let $A\in S_m(F)$ be a non-alternating
matrix. It is well-known $A$ is congruent to a non-zero diagonal
matrix, say $D$. Thus, in order to show that that the
$\GL_m(F)$-submodule of $S_m(F)$, say $M$, generated by $A$ and
$A_m(F)$ is equal to $S_m(F)$, it suffices to show that $M$
contains a matrix of rank 1. Suppose $D$ has rank $>1$; by scaling
$D$ we may assume that its first two diagonal entries are equal to
1. Let $D_1$ be the matrix obtained from $D$ by replacing its top
left $2\times 2$ corner by $\left(\begin{matrix}
0 & 1\\
1 & 1\end{matrix}\right)$ and let $D_2\in A_m(F)$ be the direct
sum of $\left(\begin{matrix}
0 & 1\\
1 & 0\end{matrix}\right)$ and the zero block. Then $D_1$ and $D_2$
belong to $M$ and $D_1-D_2$ is a non-zero diagonal matrix whose
first entry is 0. It follows by induction that $M$ contains a
matrix of rank 1, as required.

\begin{note} It might seem that {\em all} $\GL_m(F)$-submodules of
$M_m(F)$ can be obtained from above, but there is at least one
exception. If $F=F_2$ then $S_2(F)$ is the direct sum of the
1-dimensional submodule $A_m(F)$ with the 2-dimensional submodule
generated by the identity matrix.
\end{note}

\subsection{The structure of the $FG$-module $G[\VI]/(G[\VI]\cap
G[V/\R(V)])$}

We refer the reader to the definition of the $FG$-submodules
$(i)V$ of $V$, and the irreducible $FG$-modules $Q_{2k}^i$ built
upon them, both of which are defined prior to Theorem \ref{irrvI}.

\begin{notn} If $U$ and $W$ are $FG$-modules, let
$\Bil(U,W)$ denote the $FG$-module of all bilinear forms $U\times
W\to F$.
\end{notn}

\begin{defn} For $1\leq i\leq \t+1$ let
$$
\m_i=\{\phi\in\Bil(V/\VI)\,|\, (i-1)V/(0)V\subseteq \R(\phi)\}.
$$
\end{defn}

We have a series of $FG$-modules
$$
\Bil(V/\VI)=\m_1\supset \m_2\supset\cdots\supset \m_{\t+1}=(0).
$$
We further refine each link $\m_i\supset \m_{i+1}$, $1\leq i\leq
\t$ of this chain as follows.

\begin{defn}
\label{series}
 For $1\leq i\leq j\leq \t+1$ let
$$
\l_{i,j}=\{\phi\in \m_i\,|\, \phi((i)V/(0)V, (j-1)V/(0)V)=0\text{
and } \phi((j-1)V/(0)V, (i)V/(0)V)=0 \}.
$$
\end{defn}

We have a series of $FG$-modules
$$
\m_i=\l_{i,i}\supset \l_{i,i+1}\supset\cdots\supset
\l_{i,\t+1}=\m_{i+1}.
$$
This yields a refined series of $FG$-modules for $\Bil(V/\VI)$,
and intersecting each term with $S$ we get a series of
$FG$-modules for $S$. Our goal is to further refine this series
into a composition series for $S$, with known factors, as
described below.

Let $1\leq i\leq \t$. The very definition of $\m_i$ yields an
isomorphism of $FG$-modules
\begin{equation}
\label{Z1} \m_i\to \Bil(V/(i-1)V).
\end{equation}
Post-composing (\ref{Z1}) with restriction to $(i)V/(i-1)V\times
(i)V/(i-1)V$ yields a homomorphism of $FG$-modules
\begin{equation}
\label{Z2} \m_i\to \Bil((i)V/(i-1)V),
\end{equation}
whose kernel is precisely $\l_{i,i+1}$. This yields a monomorphism
of $FG$-modules
\begin{equation}
\label{Z3}\m_i\cap S/(\l_{i,i+1}\cap S)\to \Bil((i)V/(i-1)V).
\end{equation}
By restricting to
$$
Q_2^i\times Q_2^i,Q_2^i\times Q_4^i,...,Q_2^i\times Q_{2s_i}^i
$$
we get a homomorphism of $FG$-modules
\begin{equation}
\label{Z4} \m_i\cap S/(\l_{i,i+1}\cap S)\to \underset{1\leq k\leq
s_i}\oplus \Bil(Q_2^i,Q_{2k}^i).
\end{equation}
We shall show below that (\ref{Z4}) is in fact an isomorphism.

Let $1\leq i<j\leq\t$. Post-composing (\ref{Z1}) with restriction
to $(i)V/(i-1)V\times (j)V/(i-1)V$ and $(j)V/(i-1)V\times
(i)V/(i-1)V$ yields a homomorphism of $FG$-modules
\begin{equation}
\label{Z5} \l_{i,j}\to \Bil((i)V/(i-1)V,(j)V/(i-1)V)\oplus
\Bil((j)V/(i-1)V,(i)V/(i-1)V).
\end{equation}
By the very nature of $\l_{i,j}$ this yields a homomorphism of
$FG$-modules
\begin{equation}
\label{Z6} \l_{i,j}\to \Bil((i)V/(i-1)V,(j)V/(j-1)V)\oplus
\Bil((j)V/(j-1)V,(i)V/(i-1)V),
\end{equation}
whose kernel is precisely $\l_{i,j+1}$. This yields a monomorphism
of $FG$-modules
\begin{equation}
\label{Z7} \l_{i,j}\cap S/(\l_{i,j+1}\cap S)\to
\Bil((i)V/(i-1)V,(j)V/(j-1)V)\oplus \Bil((j)V/(j-1)V,(i)V/(i-1)V).
\end{equation}
By restricting to
$$
Q_2^i\times Q_2^j,Q_2^i\times Q_4^j,...,Q_2^i\times Q_{2s_j}^j
$$
and
$$
Q_2^j\times Q_2^i,Q_2^j\times Q_4^i,...,Q_2^j\times Q_{2s_i}^i.
$$
we get a homomorphism of $FG$-modules
\begin{equation}
\label{Z8} \l_{i,j}\cap S/(\l_{i,j+1}\cap S)\to \underset{1\leq
l\leq s_j}\oplus \Bil(Q_2^i,Q_{2l}^j)\bigoplus \underset{1\leq
k\leq s_i}\oplus \Bil(Q_2^j,Q_{2k}^i).
\end{equation}
We shall show below that (\ref{Z8}) is in fact an isomorphism.

\begin{thm}
\label{final} $G[\VI]/(G[\VI]\cap G[V/\R(V)])$ is an $FG$-module
of dimension
\begin{equation}
\label{equa}  \dim G[\VI]/(G[\VI]\cap G[V/\R(V)])=\dim
(V/\VI)(m_1+\cdots+m_\t),
\end{equation}
so the $FG$-module $G[{}^\infty V]$ has dimension $\dim
(V/\VI)(m_1+\cdots+m_t)$.

The irreducible constituents of $G[\VI]/(G[\VI]\cap G[V/\R(V)])$
as an $FG$-module are obtained as follows. We start with the
series for $S\cong_{FG} G[\VI]/(G[\VI]\cap G[V/\R(V)])$ produced
after Definition \ref{series} and then decompose each factor by
means of the maps (\ref{Z4}) and (\ref{Z8}), both of which are
isomorphisms.

Each summand in (\ref{Z8}) is an irreducible $FG$-module, while
the summands in (\ref{Z4}) has the constituents indicated in
Theorem \ref{symme}. More precisely, we have the following
situation.

(1) If $1\leq i\neq j\leq \t$ and $1\leq l\leq s_j$ then the
composition factor $\Bil(Q_2^i,Q_{2l}^j)$ of $G[\VI]/G[\VI]\cap
G[V/\R(V)]$ is $FG$-irreducible,
$$
G[\Bil(Q_2^i,Q_{2l}^j)]\supseteq N\rtimes \underset{k\neq i,j}\Pi
E_k,
$$
where
$$
G/(N\rtimes \underset{k\neq i,j}\Pi E_k)\cong E_i\times E_j\cong
\GL_{m_i}(F)\times \GL_{m_j}(F)
$$
acts on
$$
\Bil(Q_2^i,Q_{2l}^j)\cong M_{m_i,m_j}(F)
$$
by congruence
$$
(X,Y)\cdot A=XAY',\quad X\in \GL_{m_i}(F), A\in M_{m_i,m_j}(F),
Y\in \GL_{m_j}(F).
$$

(2) If $1\leq i\leq \t$ and $1\leq k\leq s_i$ then the factor
$\Bil(Q_2^i,Q_{2k}^i)$ of the aforementioned series
$G[\VI]/(G[\VI]\cap G[V/\R(V)])$ possesses the following
properties.
$$
G[\Bil(Q_2^i,Q_{2k}^i)]\supseteq N\rtimes \underset{l\neq i}\Pi
E_l,
$$
where
$$
G/(N\rtimes \underset{l\neq i}\Pi E_l)\cong E_i\cong \GL_{m_i}(F)
$$
acts on
$$
\Bil(Q_2^i,Q_{2k}^i)\cong M_{m_i}(F)
$$
by congruence
$$
X\cdot A=XAX',\quad X\in \GL_{m_i}(F), A\in M_{m_i}(F).
$$
The irreducible constituents of $\Bil(Q_2^i,Q_{2k}^i)$ are
therefore as indicated in Theorem \ref{symme}.
\end{thm}

\noindent{\it Proof.} We first establish the inequality
\begin{equation}
\label{dimGI} \dim G[\VI]/(G[\VI]\cap G[V/\R(V)])\geq (\dim
V/\VI)(m_1+\cdots+m_\t).
\end{equation}
Recall the $F$-linear monomorphism (\ref{h6}). We easily see that
a necessary and sufficient condition for $f\in
\mathrm{Hom}_F(V_\odd^\dag,V_\infty)$ to be in its image is that
the vectors $e^{i,p}_{2k}$, $1\leq i\leq \t$, $1\leq p\leq m_i$,
$1\leq k\leq s_i$ remain $\vf$-orthogonal under $f+1$. This yields
a linear system of
$$(s_1m_1+\cdots+s_\t m_\t)^2$$
equations in
$$
(s_1m_1+\cdots+s_\t m_\t)((s_1+1)m_1+\cdots+(s_t+1)m_t)
$$
variables. Thus
\begin{equation}
\label{dimG2} \dim G[\VI]\geq (s_1m_1+\cdots+s_\t
m_\t)(m_1+\cdots+m_t)=(\dim V/\VI)(m_1+\cdots+m_t).
\end{equation}
But from Lemma \ref{isomo} we know that
\begin{equation}
\label{dimG3} \dim G[\VI]\cap G[V/\R(V)]=\dim
\mathrm{Hom}_F(V/\VI,\R(V)).
\end{equation}
By combining (\ref{dimG2}) and (\ref{dimG3}) we obtain
(\ref{dimGI}).

We next explicitly describe the linear system governing the image
of (\ref{h6}). Let $f\in\mathrm{Hom}_F(V_\odd^\dag,V_\infty)$ and
write
\begin{equation}
\label{Xij} f(e_{2k}^{i,p})=\underset{1\leq j\leq t}\sum\;
\underset{1\leq q\leq m_i}\sum\; \underset{0\leq l\leq s_i}\sum
{}_{i,j}X_{2k,2l+1}^{p,q} e_{2l+1}^{j,q},
\end{equation}
where ${}_{i,j}X_{2k,2l+1}^{p,q}\in F$. Then
$f=(g-1)|_{V_\odd^\dag}$ for some $g\in G[\VI]$ if and only if for
$1\leq i,j\leq \t$, $1\leq k\leq s_i$, $1\leq l\leq s_j$, $1\leq
p\leq m_i$ and $1\leq q\leq m_j$ we have
\begin{equation}
\label{2Xij} 0=\langle
(f+1)(e_{2k}^{i,p}),(f+1)(e_{2l}^{j,q})\rangle={}_{i,j}X_{2k,2l+1}^{p,q}+{}_{j,i}X_{2l,2k-1}^{q,p}.
\end{equation}

We next utilize (\ref{Xij}) and (\ref{2Xij}) to show that equality
prevails in (\ref{dimGI}), and to infer from it that (\ref{Z4})
and (\ref{Z8}) are isomorphisms.

Suppose first that $1\leq i\leq \t$ and $\phi\in
\Bil((i)V/(i-1)V)$ belongs to the image of (\ref{Z3}). From the
very definition of $S$ we see that $\phi$ is the image under
(\ref{Z2}) of $\widehat{\vf_g}\in \m_i$ for some $g\in G[\VI]$.
For $1\leq k,l\leq s_i$ let ${}_{i}A_{2k,2l}\in M_{m_i}(F)$ denote
the Gram matrix of $\phi|_{Q_{2k}^i\times Q_{2l}^i}$ relative to
the bases of $Q_{2k}^i$ and $Q_{2l}^i$ described in Theorem
\ref{irrvI}. For $1\leq p,q\leq m_i$ let ${}_{i}A_{2k,2l}^{p,q}$
denote the $(p,q)$-entry of ${}_{i}A_{2k,2l}$. Then
\begin{equation}
\label{grande}
\begin{aligned}
{}_{i}A_{2k,2l}^{p,q} &=
\phi(e_{2k}^{i,p}+(i-1)V,e_{2l}^{i,q}+(i-1)V))\\
&=\widehat{\vf_g}(e_{2k}^{i,p}+\VI,e_{2l}^{i,q}+\VI)\\
&=\vf_g(e_{2k}^{i,p},e_{2l}^{i,q})\\
&=\vf((g-1)e_{2k}^{i,p},e_{2l}^{i,q})=\la
(g-1)e_{2k}^{i,p},e_{2l}^{i,q}\ra.
\end{aligned}
\end{equation}
Let $f=(g-1)|_{V_\odd^\dag}$ and let (\ref{Xij}) be the
representation of $f$ relative to our chosen basis of
$V_\odd^\dag$. Then (\ref{Xij}) and (\ref{grande}) yield
\begin{equation}
\label{3Xij} {}_{i}A_{2k,2l}^{p,q}=\langle
(g-1)e^{i,p}_{2k},e^{j,q}_{2l}\rangle= \langle
f(e^{i,p}_{2k}),e^{j,q}_{2l}\rangle= {}_{i,i}X^{p,q}_{2k,2l+1}.
\end{equation}
Applying (\ref{3Xij}) and (\ref{2Xij}), we see that, if $k>1$ then
$$
{}_{i}A_{2k,2l}^{p,q}={}_{i,i}X^{p,q}_{2k,2l+1}=-
{}_{i,i}X^{q,p}_{2l,2k-1}=-{}_{i}A_{2l,2k-2}^{q,p}.
$$
Therefore, if $k>1$ then
\begin{equation}
\label{4Xij} {}_{i}A_{2k,2l}=-\left[{}_{i}A_{2l,2k-2}\right]'.
\end{equation}
But from Theorem \ref{irrvI} we know $(i)V/(i-1)V$ is the direct
sum of its $FG$-submodules $Q_{2k}^i$, so it follows from
(\ref{4Xij}) that $\phi$ is completely determined by its
restrictions to
$$
Q_2^i\times Q_2^i,Q_2^i\times Q_4^i,...,Q_2^i\times Q_{2s_i}^i.
$$
Since (\ref{Z3}) is a monomorphism, it follows from above that
(\ref{Z4}) is also a monomorphism.

Suppose next that $1\leq i<j\leq \t$ and $(\phi_i,\phi_j)\in
\Bil((i)V/(i-1)V,(j)V/(j-1)V)\oplus \Bil((j)V/(j-1)V,(i)V/(i-1)V)$
belongs to the image of (\ref{Z7}). From the very definition of
$S$ we see that $(\phi_1,\phi_2)$ is the image under (\ref{Z6}) of
$\widehat{\vf_g}\in \m_i$ for some $g\in G[\VI]$. For $1\leq k\leq
s_i$ and $1\leq l\leq s_j$, let ${}_{i,j}A_{2k,2l}\in
M_{m_i,m_j}(F)$ and ${}_{j,i}A_{2l,2k}\in M_{m_j,m_i}(F)$ denote
the Gram matrices of $\phi_1|_{Q_{2k}^i\times Q_{2l}^j}$ and
$\phi_2|_{Q_{2l}^j\times Q_{2k}^i}$ relative to the bases of
$Q_{2k}^i$ and $Q_{2l}^i$ described in Theorem \ref{irrvI}.
Reasoning as above, we deduce that, if $k>1$ then
$$
{}_{i,j}A_{2k,2l}=-\left[{}_{j,i}A_{2l,2k-2}\right]'.
$$
As above, this implies that the pair $(\phi_1,\phi_2)$ is
completely determined by the restrictions of $\phi_1$ to
$$
Q_2^i\times Q_2^j,Q_2^i\times Q_4^j,...,Q_2^i\times Q_{2s_j}^j
$$
and restrictions of $\phi_2$ to
$$
Q_2^j\times Q_2^i,Q_2^j\times Q_4^i,...,Q_2^j\times Q_{2s_i}^i.
$$
Since (\ref{Z7}) is a monomorphism, it follows from above that
(\ref{Z8}) is also a monomorphism.

By collecting all monomorphisms (\ref{Z4}) and (\ref{Z8}), and
applying them to the series for $S$ produced after Definition
\ref{series}, we obtain the inequality
$$
\dim S \leq
 \underset{1\leq i\leq \t}\sum s_im_i^2+\underset{1\leq i\neq
j\leq \t}\sum s_jm_im_j+s_im_jm_i,
$$
that is
\begin{equation}
\label{sup} \dim S \leq (s_1m_1+\cdots+s_\t
m_\t)(m_1+\cdots+m_\t)=(\dim V/\VI)(m_1+\cdots+m_\t).
\end{equation}
By combining the inequalities (\ref{dimGI}) and (\ref{sup}) we
deduce the equality (\ref{equa}) and the fact that all maps
(\ref{Z4}) and (\ref{Z8}) are isomorphisms. The remaining
assertions of the theorem are now consequence of Theorem
\ref{irrvI}.

\subsection{Dimension of $G[{}^\infty V/V_\infty]\cap
G[V_\infty]/G[\VI]$}

Recall the $F$-vector space decomposition
$V=V_\odd\oplus(V_\even\oplus V_\ndeg)\oplus V^\dag_\odd$, and
consider a basis of $V$ formed by putting together, one after
another, bases of the 3 summands in the above decomposition. We
shall identify each element of $\GIK$ with its matrix. The Gram
matrix $A$ of $\varphi$ has the form
$$
A=\left(%
\begin{array}{ccc}
  0 & 0 & A_1 \\
  0 & A_2 & 0 \\
  A_3 & 0 & 0 \\
\end{array}%
\right).
$$
By Lemma \ref{incl} if $X\in \GIK$ then
$$
X=\left(%
\begin{array}{ccc}
  1 & Y_1 & Z \\
  0 & 1 & Y_2 \\
  0 & 0 & 1 \\
\end{array}%
\right).
$$
The equation $X'AX=A$ defining $G$ translates into
\begin{equation}
\label{tuno} Y_1'A_1+A_2Y_2=0,
\end{equation}
\begin{equation}
\label{tdos} Y_2'A_2+A_3Y_1=0,
\end{equation}
\begin{equation}
\label{ttres} Z'A_1+A_3Z+Y_2'A_2Y_2=0.
\end{equation}
By Lemma \ref{step1} the conditions for $X$ to belong to $\BB$ are
$Y_1=0$, $Y_2=0$ and (\ref{ttres}).

\begin{lem}
\label{y1y} The group $\KB$ is isomorphic to the $F$-vector space
$Y$ of all pairs $(Y_1,Y_2)$ satisfying (\ref{tuno}) and
(\ref{tdos}).
\end{lem}

\noindent{\it Proof.}  Using the above notation we define the map
$\gamma:\GIK\to Y$ given by $X\mapsto (Y_1,Y_2)$. One easily
verify that $\gamma$ is a group homomorphism with kernel $\BB$. It
remains to show that $\gamma$ is surjective. Consider the linear
map $\delta:\Hom_F(V^\dag_\odd,V_\infty)\to\End_F(V^\dag_\odd)$,
which in matrix terms is given by $Z\mapsto Z'A_1+A_3Z$. By what
we mentioned above, the kernel of $\delta$ is isomorphic to $\BB$,
which by Theorem \ref{final} has dimension $\dim
(V^\dag_\odd)\times (m_1+\cdots+m_t$). It follows that the image
of $\delta$ has dimension
$$
\dim (V^\dag_\odd)\times \dim(V_\infty)-\dim(V^\dag_\odd)\times
(m_1+\cdots+m_t)=\dim (V^\dag_\odd)\times [\dim
(V_\infty)-(m_1+\cdots+m_t)]
$$
and this equals $\dim (V^\dag_\odd)\times \dim
(V^\dag_\odd)=\dim\End_F(V^\dag_\odd)$. Thus $\delta$ is
surjective, whence $\gamma$ must be surjective as well.

\begin{prop}
\label{ufa} The dimension of the $F$-vector space $\KB$ is equal
to $\dim(V_\even\oplus V_\ndeg)\times (m_1+\cdots +m_t)$.
\end{prop}

\noindent{\it Proof.} By making use to Lemma \ref{y1y} one
verifies by direct computation that any orthogonal direct
decomposition of $V_\even\oplus V_\ndeg$ resp. $V_\odd$ yields a
corresponding direct product decomposition of $F$-vector space
$\KB$. Hence we are reduced to prove this result when both
bilinear spaces $V_\even\oplus V_\ndeg$ and $V_\odd$ are
indecomposable. Thus $V_\odd$ has a basis $e_1,...,e_{2s+1}$
relative to which the Gram matrix of $\varphi$ is equal to
$J_{2s+1}(0)$ and there are two cases to be considered.

{\sc Case I:} $V_\even=(0)$ (there is no need to assume that
$V_\ndeg$ is indecomposable).

Let $f_1,...,f_n$ be a basis of $V_\ndeg$ and let $g\in\GL(V)$.
Suppose that
$$
ge_1=e_1, ge_2=e_2+u_2+v_2,ge_3=e_3,
ge_4=e_4+u_4+v_4,...,g_{2s}=e_{2s}+u_{2s}+v_{2s},ge_{2s+1}=e_{2s+1}
$$
and
$$
gf_1=f_1+a_{1,1}e_1+a_{1,3}e_3+\cdots+a_{1,2s+1}e_{2s+1},...,
gf_n=f_n+a_{n,1}e_1+a_{n,3}e_3+\cdots+a_{n,2s+1}e_{2s+1},
$$
for some $a_{ij}\in F$, $u_{2k}\in V_\ndeg$ and $v_{2l}\in
V_\infty$.

We claim that given any choice of $a_{1,1},...,a_{n,1}$ we can
find $u_{2k}\in V_\ndeg$, $v_{2l}\in V_\infty$ and all other
$a_{ij}\in F$ such that $g\in \GIK$, and, moreover, $g\BB$ will be
unique.

It suffices to find $u_{2k}\in V_\ndeg$ and all other $a_{ij}\in
F$ so that $ge_2,...,ge_{2s}$ remain orthogonal to
$gf_1,...,gf_n$, and show that the choices for these are unique.
Indeed, the proof of Lemma \ref{y1y} explains why the $v_{2l}$
will then exist to form $g\in \GIK$, and it is clear that $g\BB$
will then be unique.

In order to find the unique $u_{2k}\in V_\ndeg$ and $a_{ij}\in F$,
$j>1$, note first that $\varphi(ge_2,gf_l)=0$ translates into
$\varphi(u_2,f_l)=-a_{l1}$, $1\leq l\leq n$. As the restriction of
$\varphi$ to $V_\ndeg$ is non-degenerate, $u_2$ exists and is
unique. Secondly $\varphi(gf_l,ge_2)=0$, translates into
$a_{l3}=-\varphi(f_l,u_2)$, so all $a_{l3}$, $1\leq l\leq n$,
exist and are unique. We may now repeat this procedure to
determine $u_4$ and then all $a_{l5}$ in a unique manner, etc.

{\sc Case II:} $V_\ndeg=(0)$ and $V_\even$ has a basis
$f_1,f_2,...,f_{2n-1},f_{2n}$ relative to which the matrix of
$\varphi$ is equal to $J_{2n}(0)$.

We first consider a family of $2n$ 1-parameter subgroups of
$\GIK$. It will be obvious from the definition that non-identity
members of different 1-parameter subgroups are linearly
independent modulo $\BB$. Our family is naturally divided into two
subfamilies, say $\gamma$ and $\delta$, each of them consisting of
$n$ 1-parameter subgroups. The $\gamma$ family consists of
$\gamma_{1,a},\gamma_{3,b},...,\gamma_{2n-1,z}\in\GIK$, where
$a,b,...,z\in F$, all of which fix
$R^\infty(V_\even)=(f_2,...,f_{2n})$ pointwise, and the $\delta$
family consists of
$\delta_{2n,a},\gamma_{2n-2,b},...,\gamma_{2,z}\in\GIK$, where
$a,b,...,z\in F$, all of which fix
$L^\infty(V_\even)=(f_1,...,f_{2n-1})$. As elements of $\GIK$ they
all fix $V_\infty=(e_1,e_3,...,e_{2s+1})$ pointwise. In the
$\gamma$ family we have
$$
\gamma_{1,a}f_1=f_1+ae_1,\gamma_{1,a}e_2=e_2-af_2,\gamma_{1,a}f_3=f_3+ae_3,
\gamma_{1,a}e_4=e_4-af_4,...
$$
$$
\gamma_{3,a}f_1=f_1,\gamma_{3,a}f_3=f_3+ae_1,\gamma_{3,a}e_2=e_2-af_4,\gamma_{3,a}f_5=f_5+ae_3,
\gamma_{3,a}e_4=e_4-af_6,...
$$
with the next $\gamma_{i,a}$ similarly defined until
$$
\gamma_{2n-1,a}f_1=f_1,...,\gamma_{2n-1,a}f_{2n-3}=f_{2n-3},\gamma_{2n-1,a}f_{2n-1}=f_{2n-1}+ae_1,
$$
and
$$
\gamma_{2n-1,a}e_2=e_2-af_{2n},\gamma_{2n-1,a}e_4=e_4,...,\gamma_{2n-1,a}e_{2s}=e_{2s}.
$$
In the $\delta$ family the first member is defined by
$$
\delta_{2n,a}f_{2n}=f_{2n}+ae_{2s+1},\gamma_{2n,a}e_{2s}=e_{2s}-af_{2n-1},
$$
and
$$
 \delta_{2n,a}f_{2n-2}=f_{2n-2}+ae_{2s-1},
\gamma_{2n,a}e_{2s-2}=e_{2s-2}-af_{2n-3},...
$$
the second member by
$$
\delta_{2n-2,a}f_{2n}=f_{2n},
\delta_{2n-2,a}f_{2n-2}=f_{2n-2}+ae_{2s+1},\gamma_{2n-2,a}e_{2s}=e_{2s}-af_{2n-3},
$$
and
$$
\delta_{2n-2,a}f_{2n-4}=f_{2n-4}+ae_{2s-1},\gamma_{2n-2,a}e_{2s-2}=e_{2s-2}-af_{2n-5},...,
$$
with the next $\delta_{i,a}$ similarly defined until
$$
\delta_{2,a}f_{2n}=f_{2n},...,\delta_{2,a}f_{4}=f_{4},
\delta_{2,a}f_{2}=f_{2}+ae_{2s+1},\delta_{2,a}e_{2s}-af_1,\delta_{2,a}e_{2s-2}=e_{2s-2},...,
\delta_{2,a}e_2=e_2.
$$

This explicit family of $2n$ 1-parameter subgroups of $\GIK$ show
that the dimension of $\KB$ is at least $2n$. We next show the
reverse inequality. For this purpose we consider the bilinear
space $W=L^2(V)/L(V)$, whose bilinear form is the one naturally
induced by $\varphi$ (this works since $L(V)$ is contained in
radical of $L^2(V)$). The canonical form-preserving linear map
$V\to W$ induces a canonical group homomorphism $G(V)\to G(W)=P$.
The latter maps $V_\infty$ into $W_\infty$ and ${}^\infty V$ into
${}^\infty W$, thereby yielding a group homomorphism, actually a
linear map from $\KB$ into $P[W_\infty]\cap P[{}^\infty
W/W_\infty]$. One verifies that the kernel of this map is
generated by the classes modulo $\BB$ of $\gamma_{2n-1,a}$ and
$\delta_{2,b}$ as $a,b$ run through $F$, so it has dimension 2.
Applying this procedure repeatedly until $\dim W_\odd=1$ or
$W_\even=0$ -in which cases our result is obvious- it follows that
$\dim \GIK\leq 2n$, as required.

As a corollary of Theorem \ref{final} and Proposition \ref{ufa} we
finally obtain

\begin{thm}
\label{ulto} $\dim\GIK=\dim (V/V_\infty)\times (m_1+\cdots+m_t)$.
\end{thm}

We know from Lemma \ref{step3} that $\GIK$ is a nilpotent group of
class $\leq 2$. The following result describes the exact
nilpotency class. The proof, which will be omitted, consists of a
case by case analysis, all of which is direct consequence of the
preceding material. We make however one clarifying remark: if
$V_\odd\neq \R(V)$ and $V_\even\neq (0)$ then the elements
$\gamma_{1,a}$ and $\delta_{2,b}$ of $\GIK$ do not commute
provided $a,b\in F$ are non-zero.

\begin{lem}
\label{cla} (a) If $V_\odd=(0)$ or $V=V_\odd$ then $\GIK$ is
trivial.

(b) If ($V_\odd\neq (0)$ and $V\neq V_\odd$) and [($V_\odd=\R(V)$)
or ($V_\even=(0)$ and $V_\odd$ has at most one indecomposable
block of size $\geq 3$ and $\dim V_\ndeg=1$)]  then $\GIK$ is
non-trivial and abelian.

(c) In all other cases $\GIK$ is non-abelian.

\end{lem}


\section{Decomposing $G(V)$ in terms of $G(V_\odd)$, $G(V_\even)$ and
$G(V_\ndeg)$}

The next result summarizes what we know about the $G(V_\odd)$. A
notable fact is that even though $V_\odd$ is far from being
uniquely determined by $V$, the image of the restriction group
homomorphism $G(V_\odd)\to GL(V_\infty)$ is the same for all
choices of $V_\odd$, as it coincides with the image of $G\to
GL(V_\infty)$.

\begin{thm}
\label{kop} We have
\begin{equation}
\label{pre} G(V_\odd)=G[\VI]\rtimes G(V_\odd)\cap
G\{V_\odd^\dag\}=G[\VI]\rtimes (U\rtimes E),
\end{equation}
where the action of $E\cong \underset{1\leq i\leq
t}\Pi\GL_{m_i}(F)$ on the unipotent group $U$, and the action of
$U\rtimes E$ on the abelian unipotent group $G[\VI]$ possess the
properties previously described in the paper.

Moreover, the restriction maps $G(V_\odd)\to \GL(V_\infty)$ and
$G\to \GL(V_\infty)$ have exactly the same image, say $H$. Indeed,
both maps restricted to $U\rtimes E$ yield the isomorphism
$U\rtimes E\to H$, while both maps have split kernels,
respectively equal to $G[\VI]$ and $G[V_\infty]$. Thus $H$ is
isomorphic to
$$
G(V_\odd)/G[\VI]\cong U\rtimes E\cong G/G[V_\infty].
$$
\end{thm}

\noindent{\it Proof.} Applying Lemma \ref{step1} to the
decomposition (\ref{26-2}) with $V=V_\odd$ and making use of
(\ref{26-1}) we get (\ref{pre}). Again by Lemma \ref{step1}, the
restriction map $G(V_\odd)\to \GL(V_\infty)$ has $G[\VI]$ in its
kernel. Let $H$ denote its image. By Lemma \ref{incl}
$G[V_\infty]\cap G(V_\odd)\cap G\{V_\odd^\dag\}=\,<1>$, whence
$U\rtimes E\to H$ is an isomorphism. It follows from (\ref{26-2})
that the image of $G\to \GL(V_\infty)$ coincides with the image of
$U\rtimes E\to \GL(V_\infty)$, that is $H$. This completes the
proof.

Next we produce further decompositions for $G$.

\begin{thm}
\label{kop2} We have the following decompositions for $G(V)$.
$$G(V)=(G[\VI/V_\infty]\cap G[V_\infty]) (G(V_\odd)\times
G(V_\even)\times G(V_\ndeg)),$$ where the intersection of
$G(V_\odd)\times G(V_\even)\times G(V_\ndeg)$ with the normal
subgroup $G[\VI/V_\infty]\cap G[V_\infty]$ of $G(V)$ is the normal
subgroup $G[\VI]$ of $G(V)$;
$$
G(V)=G[\VI/V_\infty]\cap G[V_\infty]\rtimes ((G(V_\odd)\cap
G\{V_\odd^\dag\})\times G(V_\even)\times G(V_\ndeg)),
$$
where $G(V_\odd)/G[\VI]\cong G(V_\odd)\cap G\{V_\odd^\dag\}$;
$$
G(V)=G[V_\infty] G(V_\odd),
$$
where $G[V_\infty]\cap G(V_\odd)=G[\VI]$;
$$
G=G[V_\infty]G[\VI/V_\infty],
$$
where $G[V_\infty]\cap G[\VI/V_\infty]$ is a unipotent normal
subgroup of $G$ with nilpotency class $\leq 2$.
\end{thm}

\noindent{\it Proof.} The first three decompositions follow from
Theorems \ref{bigsplit} and \ref{kop}, while the fourth follows
from the third.

Finally we consider the special but interesting case when
$V=V_\odd$ is {\em homogenous}, namely when $V=V_{\odd}$ is the
direct sum of $m$ Gabriel blocks of equal size $2s+1$. The
isomorphism type of $G$ is fully revealed in this case.



\begin{thm}
\label{homo} Suppose $V=V_{\odd}$ is the direct sum of $m$ Gabriel
blocks of size $2s+1$. Then
$$G\cong (\underset{1\leq k\leq s}\Pi M_m(F))\rtimes
\GL_m(F),$$ where $\GL_m(F)$ acts diagonally on $\underset{1\leq
k\leq s}\Pi M_m(F)$ by congruence.

Internally, $G[V_\infty]$ has a natural structure of $FG$-module
of dimension $sm^2$ over $F$. As a module over $\GL_m(F)\cong
G/G[V_\infty]$, $G[V_\infty]$ is isomorphic to $\underset{1\leq
k\leq s}\Pi M_m(F)$, upon which $\GL_m(F)$ acts diagonally by
congruence.
\end{thm}

\noindent{\it Proof.} Observe first of all that
$N=G[V_\infty]=G[\VI]$, so $G=G[\VI]\rtimes E$, where
$E\cong\GL_{m}(F)$ and the action of $E$ on $G[\VI]$ has been
determined. More precisely, as the case $s=0$ is obvious, we may
assume that $s\geq 1$. Since $t=1$ and $\R(V)=0$, Theorem
\ref{final} yields
$$
G[\VI]=G[\VI]\cap G[V/\R(V)]\cong S\cong\underset{1\leq k\leq
s}\oplus\Bil(Q_2^1,Q_{2k}^1),
$$
and the indicated action of $E\cong M_m(F)$ on $\underset{1\leq
k\leq s}\oplus\Bil(Q_2^1,Q_{2k}^1)\cong \underset{1\leq k\leq
s}\Pi M_m(F).$

\section{Inductive approach}

It is possible to extract useful information on $G$ by studying
the canonical group homomorphism
$$G(V)\to G(L^2(V)/L^1(V)).$$
Here the bilinear space $L^2(V)/L^1(V)$ can be obtained from $V$
in a straightforward manner: its non-degenerate parts are
equivalent, and all Gabriel blocks of $V$ decrease in size by 2
when passing from $V$ to $L^2(V)/L^1(V)$, except for those of size
$\leq 2$ which disappear. The above map is very likely to be
surjective (we have checked this in a few cases), so repeated
application of it would yield $G$ as constructed from $G(V_\ndeg)$
and the various kernels, all of which respond to the same pattern.

This sort of approach seems to be applicable to $G(V_\ndeg)$, once
it is already decomposed as in (\ref{dw}). There is a canonical
$G$-invariant filtration for $V_\ndeg$ and one produces from it a
non-degenerate bilinear space as a section of $V_\ndeg$. Special
cases have revealed the associated group homomorphism to be
surjective as well.





\end{document}